\documentclass[11pt,a4paper]{article}

\usepackage{amsmath,amssymb}
\usepackage{color}
\newcommand{\C}{\mathbb{C}}
\newcommand{\R}{\mathbb{R}}
\newcommand{\N}{\mathbb{N}}
\newcommand{\Z}{\mathbb{Z}}
\def\cA{{\mathcal A}}
\def\cC{{\mathcal C}}

\def\cS{{\mathcal S}}
\newcommand{\ee}{\varepsilon}
\renewcommand{\aa}{\alpha}
\renewcommand{\div}{{\rm div}\,}
\newcommand{\curl}{{\rm curl}\,}

\newcommand{\Frac}{\displaystyle \frac}

\newcommand{\Sum}{\displaystyle \sum}

\def\d{\partial}
\def\ddl{\dot \Delta_l}
\def\ddj{\dot \Delta_j}
\def\ddq{\dot \Delta_q}

\def\tilde{\widetilde}
\def\hat{\widehat}

\newcommand{\fd}{\frac{d}{2}}

\newcommand{\ue}{u_\ee}
\newcommand{\qa}{q_\aa}
\newcommand{\ua}{u_\aa}
\newcommand{\ca}{c_\aa}
\newcommand{\phia}{\phi_\aa}
\newcommand{\dq}{\delta q}
\newcommand{\du}{\delta u}

\newtheorem{thm}{Theorem}
\newtheorem{lem}{Lemma}

\newtheorem{prop}{Proposition}
\newtheorem{defi}{Definition}

\newtheorem{rem}{Remark}

\title{Convergence of a low order non-local Navier-Stokes-Korteweg system: the order-parameter model}

\author{Fr\'ed\'eric Charve\footnote{Universit\'e Paris-Est Cr\'eteil, Laboratoire d'Analyse et de Math\'ematiques Appliqu\'ees (UMR 8050), 61 Avenue du G\'en\'eral de Gaulle, 94 010 Cr\'eteil Cedex (France). E-mail: frederic.charve@u-pec.fr}}

\date{}
\begin{document}

\maketitle

\begin{abstract} In the present article we consider a capillary compressible system introduced by C. Rohde after works of Bandon, Lin and Rogers, called the order-parameter model, and whose aim is to reduce the numerical difficulties that one encounters in the case of the classical local Korteweg system (involving derivatives of order three) or the non-local system (also introduced by Rohde after works of Van der Waals, and which involves a convolution operator).
We prove that this system has a unique global solution for initial data close to an equilibrium and we precisely study the convergence of this solution towards the local Korteweg model.
\end{abstract}
\section{Introduction}
\subsection{Presentation of the systems}
In the mathematical study of liquid-vapour mixture, Gibbs first modelled phase transitions thanks to the minimization of an energy functional with a nonconvex energy density (see \cite{Gibbs}). The phases are separated by an hypersurface and there are mainly two ways to describe it: either we consider that the interface behaves like a discontinuity for the fluid parameters (this is the Sharp Interface model), either we consider that between the phases lies a thin region of continuous transition (this is the Diffuse Interface approach, where the phase changes are seen through the variations of the density and which is much simpler numerically). Unfortunately the basic models provide an infinite number of solutions (few of them being physically relevant) and this is why authors tried to penalize the high variations of the density (with capillary terms related to surface tension) in order to select the physically correct solutions.

In the present paper, we are interested in the local and non-local Korteweg systems (in the diffuse interface model). These systems are based upon the compressible Navier-Stokes system with a Van der Waals state law for ideal fluids, and endowed with a capillary tensor.

Let us recall that the local model was introduced by Korteweg and renewed by Dunn and Serrin (see \cite{3DS}) and the non-local model was introduced by Van der Waals and renewed by F. Coquel, D. Diehl, C. Merkle and C. Rohde. For an in-depth presentation of the capillary models, we refer to \cite{Rohdehdr} and \cite{5CR}).

Let $\rho$ and $u$ denote the density and the velocity of a compressible viscous fluid ($\rho$ is a non-negative function and $u$ is a vector-valued function defined on $\R^d$). We denote by $\cA$ the following diffusion operator
$$
\cA u= \mu \Delta u+ (\lambda+\mu)\nabla \div u, \quad \mbox{with} \quad \mu>0 \quad \mbox{and} \quad \nu=\lambda+ 2 \mu >0.
$$
The Navier-Stokes equations for compressible fluids endowed with internal capillarity read:
$$
\begin{cases}
\begin{aligned}
&\d_t\rho+\div (\rho u)=0,\\
&\d_t (\rho u)+\div (\rho u\otimes u)-\cA u+\nabla(P(\rho))=\kappa\rho\nabla D[\rho].\\
\end{aligned}
\end{cases}
$$
The capillary coefficient $\kappa$ may depend on $\rho$ but in this article it is chosen constant. In the local Korteweg system $(NSK)$, the capillary term $D[\rho]$ is given by (see \cite{3DS}):
$$D[\rho]=\Delta\rho,$$
and, in the non-local Korteweg system $(NSRW)$ (see \cite{5Ro}, \cite{5CR}, and \cite{VW}), if $\phi$ is an interaction potential which satisfies the following conditions
\begin{equation}
 (|.|+|.|^2)\phi(.)\in L^1(\R^d)\mbox{, }\quad\int_{\R^d}\phi(x)dx=1,\quad\phi\mbox{ even, and }\phi\geq0,
\label{potentielcond}
\end{equation}
then $D[\rho]$ is the non-local term given by:
$$D[\rho]=\phi*\rho-\rho.$$

Comparing the Fourier transform of the capillary terms, we have $(\hat{\phi}(\xi)-1) \hat{\rho}(\xi)$ in the non-local model, and $-|\xi|^2 \hat{\rho}(\xi)$ in the local model so that a natural question is to study the closedness of the solutions of these models when $\hat{\phi}(\xi)$ is formally "close" to $1-|\xi|^2$. For this, we introduced in \cite{CH} a specific interaction potential and considered the following non-local system: 
$$
\begin{cases}
\begin{aligned}
&\d_t\rho_\ee+\div (\rho_\ee \ue)=0,\\
&\d_t (\rho_\ee \ue)+\div (\rho_\ee u\otimes \ue)-\cA \ue+\nabla(P(\rho_\ee))=\rho_\ee\frac{\kappa}{\ee^2}\nabla(\phi_\ee*\rho_\ee-\rho_\ee),\\
\end{aligned}
\end{cases}
\leqno{(NSRW_\ee)}
$$
with
$$\phi_{\ee}=\frac{1}{\ee^d} \phi(\frac{x}{\ee}) \quad \mbox{with} \quad \phi(x)=\frac{1}{(2 \pi)^d} e^{-\frac{|x|^2}{4}}$$
For a fixed $\xi$ the Fourier transform of $\phi_\ee$ is $\hat{\phi_\ee}(\xi)=e^{-\ee^2|\xi|^2}$, and when $\ee$ is small, $\Frac{\hat{\phi_\ee}(\xi)-1}{\ee^2}$ is close to $-|\xi|^2$.

Using energy methods, we proved in \cite{CH} that this system has a unique global strong solution for initial data close to an equilibrium state. The functional setting are classical and hybrid Besov spaces (taylored to the capillary term). We also obtained that when the small parameter $\ee$ goes to zero, the solution tends to the corresponding solution of the local Korteweg system and we obtained a rate of convergence in terms of $\ee$. In \cite{CHVP}, we provided by Lagrangian methods more precise a priori estimates giving a better understanding of the convergence and the hybrid Besov setting in terms of the linear Fourier structures.

Though, these models are not completely satisfying. On one hand, recalling the results from \cite{Dinv}, \cite{Dbook} (compressible Navier-Stokes system), \cite{DD}, \cite{CH1} (local Korteweg model), and \cite{Has1}, \cite{CH} (non-local Korteweg model), we observe that the density in the local capillary model is far more regular than in the non-local model where it shares the same frequency structure as in the compressible Navier-Stokes model (heat regularization in low frequencies and only a damping in high frequencies).

On the other hand, from a numerical point of view the local model is difficult to handle because the capillary term contains third-order derivatives. The non-local model also presents difficulties in numerical studies: even if the capillary term only contains derivatives of order one, it involves a convolution operator, whose numerical difficulty is comparable.

For this reason, C.Rohde presented in \cite{Rohdeorder} a new model, called the order-parameter model, and inspired by the work of D. Brandon, T. Lin and R. C. Rogers in \cite{BLR}.

This new system consists in introducing in the capillary term $\alpha^2 \nabla (c-\rho)$ a new variable $c$ called the "order parameter", which is coupled to the density via the following relation related to the Euler-Lagrange equation from the variational approach ($\aa$ controls the coupling between $\rho$ and $c$):
$$
\ee^2 \Delta c +\alpha^2 (\rho-c)=0
$$
so that the new system he considers is the following:
$$
\begin{cases}
\begin{aligned}
&\d_t\rho_\aa+\div (\rho_\aa \ua)=0,\\
&\d_t (\rho_\aa \ua)+\div (\rho_\aa \ua\otimes \ua)-\cA \ua+\nabla(P(\rho_\aa))=\kappa \aa^2 \rho_\aa\nabla(\ca-\rho_\aa),\\
&\ee^2 \Delta \ca +\alpha^2 (\rho_\aa -\ca)=0.
\end{aligned}
\end{cases}
\leqno{(NSOP_\aa)}
$$
As emphasized by C. Rohde, from a numerical point of view this system is much more interesting because now we only have one derivative in the capillary tensor (which is local), and the additionnal equation for the order parameter is a simple linear elliptic equation that can be easily and numerically fast solved at least when the mesh is fixed. Moreover as we will see later, as for the previous non-local capillary model, this system has the same frequency structure as the classical Navier-Stokes model.

In \cite{Rohdeorder}, C. Rohde proves (for $\ee=\lambda=\mu=1$ in the two-dimensionnal case) that the system has a unique local classical solution:

\begin{thm}[\cite{Rohdeorder}]
\sl{Assume that the initial data $(\rho_0, u_0)$ is independant of $\aa>0$ with $u_0\in H^4(\R^2)$, $\rho_0>0$ and $\rho_0-\bar{\rho}\in H^4(\R^2)$ for some constant $\bar{\rho}$. Let $c_0$ the solution of the elliptic problem $-\Delta c_0+ \aa^2 c_0=\aa^2 \rho_0$. There exists a constant $T_*>0$ such that the initial-value problem $(NSOP_\aa)$ has a unique solution defined on $[0, T_*[$ satisfying:
$$\rho_\aa-\bar{\rho}, \ua \in L^\infty(0,T_*; H^4(\R^2)), \quad \rho_\aa>0,$$
$$\ca-\bar{\rho} \in L^\infty(0,T_*; H^5(\R^2)).$$
Moreover, for all $t\in[0,T_*[$ we have
$$
\lim_{\aa \rightarrow \infty} \|\rho_\aa (t,.)-\ca(t,.)\|_{L^2(\R^2)}=0.
$$
}
\end{thm}
In this paper, C. Rohde also conjectures that when the coupling constant $\aa$ goes to infinity, these solutions converge to the solution of the local Korteweg model.

\subsection{Statement of the results}

In the present article, following what we did in the whole space for the non-local system (see \cite{CH}) and using lagrangian methods from \cite{CHVP}, we will prove that under smallness conditions, and with less regular initial data, the system has global strong solutions in the following critical spaces (we refer to the appendix for more details on Besov spaces and hybrid spaces). We also prove the above conjectured convergence and give an explicit rate of convergence with respect to $\alpha$.
\begin{defi}
\sl{The space $F_{\aa}^s$ is the set of functions $(q,c,u)$ in
$$
\left(\cC_b(\R_+, \dot{B}_{2,1}^{s-1}\cap \dot{B}_{2,1}^s)\cap L^1(\R_+, \dot{B}_\aa^{s+1,s-1}\cap \dot{B}_\aa^{s+2,s})\right)^2 \times
\left(\cC_b(\R_+, \dot{B}_{2,1}^{s-1})\cap L^1(\R_+, \dot{B}_{2,1}^{s+1})\right)^d
$$
endowed with the norm $\|(q,c,u)\|_{F_\aa^s}=\|(q,c,u)\|_{F_\aa^s(\infty)}$ where for all $t$ we denote (recall that $\nu_0=\min(\mu, \nu)$)
\begin{multline}
\|(q,c,u)\|_{F_\aa^s(t)} \overset{def}{=}  \|u\|_{\tilde{L}_t^{\infty} \dot{B}_{2,1}^{s-1}}+ \|q\|_{\tilde{L}_t^{\infty} \dot{B}_{2,1}^{s-1}}+ \nu\|q\|_{\tilde{L}_t^{\infty} \dot{B}_{2,1}^{s}} +\|c\|_{\tilde{L}_t^{\infty} \dot{B}_{2,1}^{s-1}}+ \nu\|c\|_{\tilde{L}_t^{\infty} \dot{B}_{2,1}^{s}}\\
+\nu_0\|u\|_{\tilde{L}_t^1 \dot{B}_{2,1}^{s+1}}+ \nu\|q\|_{\tilde{L}_t^1 \dot{B}_{\aa}^{s+1,s-1}}+ \nu^2\|q\|_{\tilde{L}_t^1 \dot{B}_{\aa}^{s+2,s}} +\nu\|c\|_{\tilde{L}_t^1 \dot{B}_{\aa}^{s+1,s-1}}+ \nu^2\|c\|_{\tilde{L}_t^1 \dot{B}_{\aa}^{s+2,s}}
\end{multline}}
\end{defi}
\begin{thm}
\sl{Let $\aa>0$ and assume $\min(\mu,2\mu+\lambda)>0$. There exist two positive constants $\eta_{OP}$ and $C$ only depending on $d$, $\mu$, $\lambda$, $\kappa$, and $P'(\overline{\rho})$ such that for all $\eta\leq \eta_{OP}$, if $\rho_0-\overline{\rho}\in \dot{B}_{2,1}^{\fd-1}\cap \dot{B}_{2,1}^{\fd}$, $u_0 \in \dot{B}_{2,1}^{\fd-1}$, $c_0$ is defined by $-\Delta c_0+ \aa^2 c_0=\aa^2 \rho_0$ and
$$
\|\rho_0-\overline{\rho}\|_{\dot{B}_{2,1}^{\fd-1} \cap \dot{B}_{2,1}^{\fd}} +\|u_0\|_{\dot{B}_{2,1}^{\fd-1}}\leq \eta
$$
then system $(NSOP_\aa)$ has a unique global solution $(\rho_\aa,\ca, \ua)$ with $(\rho_\aa-\overline{\rho}, \ca-\overline{\rho},\ua)\in F_{\aa}^{\fd}$ such that:
$$
\|(\rho_\aa-\overline{\rho},\ca-\overline{\rho},\ua)\|_{F_{\aa}^\fd} \leq C^0 \overset{def}{=} C (\|\rho_0-\overline{\rho}\|_{\dot{B}_{2,1}^{\fd-1} \cap \dot{B}_{2,1}^{\fd}} +\|u_0\|_{\dot{B}_{2,1}^{\fd-1}}).
$$
Moreover we have the global in time results:
$$
\begin{cases}
\vspace{0.2mm}
\|\ca-\rho_\aa\|_{\tilde{L}^\infty(\R_+, \dot{B}_{2,1}^{\fd-1})} +\nu \|\ca-\rho_\aa\|_{\tilde{L}^\infty(\R_+,  \dot{B}_{2,1}^{\fd})} \underset{\aa\rightarrow \infty}{\longrightarrow} 0,\\
\nu \|\ca-\rho_\aa\|_{L^1(\R_+, \dot{B}_{2,1}^{\fd-1})} +\nu^2 \|\ca-\rho_\aa\|_{L^1(\R_+, \dot{B}_{2,1}^{\fd})} \leq C^0\aa^{-2}.
\end{cases}
$$
}
\label{thexist}
\end{thm}
The following result deals with the convergence in $\aa$: when the initial data are small enough (so that we have global solutions for $(NSK)$ and $(NSOP_\aa)$) the solution of $(NSOP_\aa)$ goes to the solution of $(NSK)$ when $\aa$ goes to infinity.

\begin{thm}
\sl{With the same assumptions as before, there exists $0<\eta_0\leq \min(\eta_K, \eta_{OP})$ such that for all $\eta\leq \eta_0$, if
$$
\|\rho_0-\overline{\rho}\|_{\dot{B}_{2,1}^{\fd-1} \cap \dot{B}_{2,1}^{\fd}} +\|u_0\|_{\dot{B}_{2,1}^{\fd-1}}\leq \eta,
$$
then systems $(NSK)$ and $(NSOP_\aa)$ both have global solutions and $\|(\rho_\aa-\rho, \ca-\rho, \ua-u)\|_{F_\aa^{\fd}}$ goes to zero as $\aa$ goes to infinity. Moreover, with the same notations as before, there exists a constant $C=C(\eta, \kappa, \overline{\rho}, P'(1))>0$ such that for all $h\in ]0, 1[$ (if $d=2$) or $h\in ]0,1]$ (if $d\geq 3$)
$$
 \|(\rho_\aa-\rho, \ca-\rho, \ua-u)\|_{F_\aa^{\fd-h}}\leq C \aa^{-h},
$$
}
\label{thcv}
\end{thm}
\begin{rem}
\sl{We can assume that $\ee=1$ without loss of generality. If not we just have to replace $\kappa$ by $\ee^2$ and $\aa$ by $\aa/\ee$.
}
\end{rem}
\begin{rem}
\sl{
As the order parameter $\ca$ goes to $\rho_\aa$, we formally get that when $\aa$ goes to infinity, the capillary term goes to $\kappa \rho \nabla \Delta \rho$.
}
\end{rem}

\subsection{outline of the paper}

The article is structured the following way: section $2$ is devoted to the proof of theorem \ref{thexist}. We first introduce an interaction potential $\phia$ that allows us to rewrite the system into a non-local shape. As we want precise estimates we follow the methods from \cite{CHVP}: we first obtain estimates on the linearized system and then on the advected linear system thanks to a Lagrangian change of variable. The rest of the proof is classical, we define approximated solutions thanks to the Friedrichs' scheme and obtain existence and uniqueness like in \cite{CH}. In section $3$ we prove theorem \ref{thcv} and in the appendix, we first recall basic properties of Besov spaces, then we provide estimates for the flow of a smooth vectorfield. The last part of the appendix is devoted to Bessel functions that are needed for the expression of our new interaction potential.

\section{Proof of theorem \ref{thexist}}

\subsection{Interaction potential}

As announced in the introduction, we first rewrite the system in a non-local shape. Let us focus on the last equation, we can write that (for more clarity we drop the subscripts with $\aa$):
$$
-\ee^2 \Delta (\rho-c) +\alpha^2 (\rho -c)=-\ee^2 \Delta \rho
$$
which leads to:
$$
\aa^2 (\rho-c)=-\aa^2 (-\Delta +\frac{\aa^2}{\ee^2})^{-1} \Delta \rho=- (\frac{-\Delta}{\aa^2} +\frac{1}{\ee^2})^{-1} \Delta \rho
$$
So that in Fourier variable:
$$
\widehat{\aa^2(\rho-c)}(\xi)=\frac{|\xi|^2}{\frac{|\xi|^2}{\aa^2}+\frac{1}{\ee^2}} \hat{\rho}(\xi)= \ee^2 \cdot \frac{\aa^2}{\ee^2}(1-\frac{1}{\frac{\ee^2}{\aa^2} |\xi|^2+1})\hat{\rho}(\xi).
$$
Then up to choose $\kappa=\ee^2$ and replace $\aa$ by $\aa/\ee$, from now on we assume that $\ee=1$ and then if we introduce:
$$
D[\rho]=\aa^2 (c-\rho),
$$
we have
\begin{equation}
\widehat{D[\rho]}(\xi)=\frac{-|\xi|^2}{\frac{|\xi|^2}{\aa^2}+1} \hat{\rho}(\xi) =\aa^2 (\frac{1}{\frac{|\xi|^2}{\aa^2}+1}-1)\hat{\rho}(\xi).
\label{capillorder}
\end{equation}
As a consequence, when $\aa$ is large, $\aa^2(\rho-c)$ formally goes to $\Delta \rho$ as for the non-local capillary term from \cite{CHVP}, and the object of this article is to prove that the solutions of this system will go to the solutions of the local Korteweg model. Let us now define the interaction potential $\phia$ by:
\begin{equation}
\widehat{\phia}(\xi)= \frac{1}{\frac{|\xi|^2}{\aa^2}+1}.
\label{phia}
\end{equation}
We have $\int_\R \phia(x)dx=1$ and $D[\rho]=\aa^2(\phia *\rho -\rho)$. If we put $\phi=\phi_1$ then
\begin{equation}
\widehat{\phi}(\xi)=\frac{1}{|\xi|^2+1}, \quad \widehat{\phia}=\widehat{\phi}(\cdot /\aa), \quad \mbox{and }\phia=\aa^d \phi(\aa \cdot).
\end{equation}
In some cases we have explicit expressions for this inverse Fourier transform: for all $x$, $\phi(x)=C e^{-|x|}$ when $d=1$, $\phi(x)=C' \frac{e^{-|x|}}{|x|}$ when $d=3$ (we refer to \cite{Stein}). In the other cases the expression of $\phi$ involves Bessel functions. Let us begin by recalling that the fourier transform of a radial function is also radial, more precisely (see for example \cite{Stein} page 213) there exists a constant $C_d$ such that if $f(x)=f_0 (|x|)$ for all $x\in \R^d$, then its Fourier transform satisfies for all $\xi \in\R^d$, $\hat{f}(\xi)=F_0(|\xi|)$ where for all $\rho>0$
$$
F_0(\rho)=\frac{C_d}{\rho^{\fd-1}} \int_0^\infty J_{\fd -1}(\rho r) f_0(r) r^{\fd} dr,
$$
where $J_\nu$ denotes the general Bessel function of real index $\nu$. This formulation is related to the Hankel transform, we refer to the appendix for more details and properties on Bessel functions.
Coming back to our problem, we then obtain that for all $x\in \R^d$,
\begin{equation}
\phi(x)=\frac{C_d}{|x|^{\fd-1}} \int_0^\infty J_{\fd -1}(r|x|) \frac{r^{\fd}}{1+r^2} dr.
\label{defphi}
\end{equation}
And thanks to the Hankel-Nicholson integrals (we refer for example to \cite{Luke} page 330 or \cite{Watson} page 434), under the following assumptions:
$$
a>0, \quad Re(z)>0, \quad -1<Re(\nu)< 2 Re(\mu)+\frac{3}{2},
$$
we have the identity:
$$
\int_0^\infty \frac{t^{\nu+1} J_\nu (at)}{(t^2+z^2)^{\mu+1}} dt =\frac{a^\mu z^{\nu-\mu}}{2^\mu \Gamma(\mu+1)} K_{\nu-\mu}(az),
$$
where $K_\nu$ denotes the modified Bessel function of the second kind and index $\nu$ (also called Hankel, Schl\"afti or Weber function). This allows us to finally write that for all $x\in \R^d$ provided that $d\in\{1,2,3,4\}$ (from the previous conditions with $\nu=\fd-1$, $\mu=0$, $z=1$, $a=|x|$),
\begin{equation}
\phi(x)=\frac{C_d}{|x|^{\fd-1}} K_{\fd-1}(|x|).
\label{phiexpl}
\end{equation}
\begin{rem}
\sl{Another way to understand the limitation on the dimension consists in observing in the integral \eqref{defphi}, that if we roughly approximate the Bessel function by $\cos(r) r^{-1/2}$ at infinity, then the integrated function has the following asymptotic expansion at infinity: $\cos(r) r^{d/2-5/2}$ ($|x|=1$ for more simplicity).
}
\end{rem}
In fact \eqref{phiexpl} is also valid for dimensions $d\geq 5$ (Like the Fourier transform, the Hankel transform can be generalized for tempered distributions). Let us compute the Fourier transform: for all $\xi\in \R^d$
$$
\int_{\R^d} e^{-ix\cdot \xi} \frac{K_{\fd-1}(|x|)}{|x|^{\fd-1}} dx =\int_0^\infty r^\fd K_{\fd-1}(r) \left( \int_{\mathbb{S}^{d-1}} e^{-ir\omega \cdot \xi} d\omega\right) dr,
$$
and thanks to the radial symmetry:
$$
\int_{\mathbb{S}^{d-1}} e^{-ir\omega \cdot \xi} d\omega = \int_{\mathbb{S}^{d-1}} e^{-ir|\xi| \omega \cdot e_1} d\omega.
$$
Performing a d-dimensional spherical change of variable we obtain that (with $\lambda=r|\xi|$):
\begin{multline}
\int_{\mathbb{S}^{d-1}} e^{-i \lambda \omega \cdot e_1} d\omega = \int_0^\pi \int_0^{2\pi}...\int_0^{2\pi} e^{-i\lambda \cos \theta_1} \sin^{d-2} \theta_1 \sin^{d-3} \theta_2 ... \sin\theta_{d-2} d\theta_1...d\theta_{d-1}\\
=C_d \int_0^\pi e^{-i\lambda \cos \theta_1} \sin^{d-2} \theta_1 d\theta_1.
\end{multline}
Using the following integral representation of function $I_\nu$ (for $Re(\nu)>-1/2$ see the appendix for modified Bessel functions $I_\nu$ and $K_\nu$)
$$
I_{\nu}(z)=\frac{z^\nu}{2^\nu \pi^{\frac{1}{2}}\Gamma(\nu+\frac{1}{2})} \int_0^\pi e^{-z\cos t} \sin^{d-2} t dt.
$$
Then thanks to the following identity (here $a=i|\xi|$ and $b=1$):
$$
\int z I_\nu(az) K_\nu(bz) dz=\frac{z}{a^2-b^2}\left(a I_{\nu+1}(az) K_\nu(bz) +bI_\nu (az) K_{\nu+1}(bz)\right),
$$
we obtain that:
$$
\int_{\R^d} e^{-ix\cdot \xi} \frac{K_{\fd-1}(|x|)}{|x|^{\fd-1}} dx =\frac{2^{\fd-1} \Gamma(\fd)}{1+|\xi|^2}.
$$
so that in \eqref{phiexpl}, $C_d=(2^{\fd-1} \Gamma(\fd))^{-1}$,
Considering the asymptotics of function $K_{\fd-1}$, $\phi$ is continuous on $\R^s-\{0\}$. Near $0$, and for $d\geq 3$, we have $\phi(x)\sim C_d |x|^{2-d}$ so that $|x|\phi(x)$ is a $L^1$ function on $\R^d$.
\label{interactionpotentielgen}
\subsection{Reformulation of the system}
We are now able to write the system into a non-local form:
$$
\begin{cases}
\begin{aligned}
&\d_t\rho_\aa+\div (\rho_\aa \ua)=0,\\
&\d_t (\rho_\aa \ua)+\div (\rho_\aa u\otimes \ua)-\cA \ua+\nabla(P(\rho_\aa))=\rho_\aa \kappa \aa^2 \nabla(\phi_\aa*\rho_\aa-\rho_\aa),\\
\end{aligned}
\end{cases}
\leqno{(NSOP_\aa)}
$$
with
$$\phia=\aa^d \phi(\aa \cdot) \quad \mbox{with} \quad \phi(x)=\frac{C_d}{|x|^{\fd-1}} K_{\fd-1}(|x|).$$
\begin{rem}
\sl{From the previous computations, we immediately get that
$$
\ca-\rho_\aa=(-\Delta+ \aa I_d)^{-1}\Delta \rho_\aa=\phia * \rho_\aa - \rho_\aa
$$
that is $\ca=\phia * \rho_\aa$. This is why we cannot choose any initial data for the order parameter and take $c_0=\phia *\rho_0$.
\label{orderparamexpr}
}
\end{rem}

As we consider initial data close to an equilibrium state $(\overline{\rho},0)$ we begin with the classical change of function $\rho= \overline{\rho}(1+q)$.  For simplicity we take $\overline{\rho}=1$. The previous system becomes (also denoted by $(NSOP_\aa)$):
$$
\begin{cases}
\begin{aligned}
&\d_t \qa+ \ua.\nabla \qa+ (1+\qa)\div \ua=0,\\
&\d_t \ua+ \ua.\nabla \ua -\cA \ua+P'(1).\nabla \qa-\kappa \aa^2 \nabla(\phia*\qa-\qa)=K(\qa).\nabla \qa- I(\qa) \cA \ua,\\
\end{aligned}
\end{cases}
\leqno{(NSOP_\aa)}
\label{NSOP}
$$
where $K$ and $I$ are the real-valued functions defined on $\R$ given by:
$$
K(q)=\left(P'(1)-\frac{P'(1+q)}{1+q}\right) \quad \mbox{and} \quad I(q)=\frac{q}{q+1}.
$$
The functional spaces we will really use are the following:
\begin{defi}
\sl{The space $E_{\aa}^s$ is the set of functions $(q,u)$ in
$$
\left(\cC_b(\R_+, \dot{B}_{2,1}^{s-1}\cap \dot{B}_{2,1}^s)\cap L^1(\R_+, \dot{B}_\aa^{s+1,s-1}\cap \dot{B}_\aa^{s+2,s})\right) \times
\left(\cC_b(\R_+, \dot{B}_{2,1}^{s-1})\cap L^1(\R_+, \dot{B}_{2,1}^{s+1})\right)^d
$$
endowed with the norm $\|(q,u)\|_{E_\aa^s}=\|(q,u)\|_{E_\aa^s(\infty)}$ where for all $t$ we denote (recall that $\nu_0=\min(\mu, \nu)$)
\begin{multline}
\|(q,u)\|_{E_\aa^s(t)} \overset{def}{=}  \|u\|_{\tilde{L}_t^{\infty} \dot{B}_{2,1}^{s-1}}+ \|q\|_{\tilde{L}_t^{\infty} \dot{B}_{2,1}^{s-1}}+ \nu\|q\|_{\tilde{L}_t^{\infty} \dot{B}_{2,1}^{s}}\\
+ \nu_0\|u\|_{\tilde{L}_t^1 \dot{B}_{2,1}^{s+1}}+ \nu\|q\|_{\tilde{L}_t^1 \dot{B}_{\aa}^{s+1,s-1}}+ \nu^2\|q\|_{\tilde{L}_t^1 \dot{B}_{\aa}^{s+2,s}}
\label{normeE}
\end{multline}}
\end{defi}
\begin{rem}
\sl{Due to obvious simplifications we slightly changed the notations for $E_{\aa}^s$ and $\dot{B}_\aa^{s+2,s}$: with the notations from \cite{CHVP} these spaces would have been respectively denoted by $E_{1/\aa}^s$ and $\dot{B}_{1/\aa}^{s+2,s}$.
}
\end{rem}
We will now follow the tracks of \cite{CH} and \cite{CHVP} to prove the results. Classically in the study in critical spaces of compressible Navier-Stokes-type systems (see \cite{Dinv, CD, arma}), the proofs of theorems \ref{thexist} and \ref{thcv} (see \cite{CH} section 2) rely on key a priori estimates on the following advected linear system ($\aa>0$ is fixed and for more simplicity we write $(q,u)$ instead of $(\qa,\ua)$):
$$
\begin{cases}
\begin{aligned}
&\d_t q+ v.\nabla q+ \div u= F,\\
&\d_t u+ v.\nabla u -\cA u+ p\nabla q-\kappa \aa^2\nabla(\phia*q-q)= G.\\
\end{aligned}
\end{cases}
\leqno{(LOP_\aa)}
$$
With
$$\cA u= \mu \Delta u+ (\lambda+\mu)\nabla \div u.$$
Although the potential function is different from the gaussian from \cite{CH}, we can easily adapt the energy methods and results from this paper. Here we will directly focus on more refined estimates as in \cite{CHVP} and use them in the proof of the last theorem: we can prove that the estimates are similar up to slight changes in the constants:
\begin{thm} Let $\aa>0$, $-\fd+1<s<\fd+1$, $I=[0,T[$ or $[0, +\infty[$ and $v\in L^1(I,\dot{B}_{2,1}^{\fd+1}) \cap L^2 (I,\dot{B}_{2,1}^{\fd})$. Assume that $(q,u)$ is a solution of System $(LOP_\aa)$ defined on $I$. There exists $\aa_0>0$, a constant $C>0$ depending on $d$, $s$ such that if $\aa\geq \aa_0$, for all $t\in I$ (denoting $\nu=\mu+2\lambda$ and $\nu_0=\min(\nu, \mu)$),
\begin{multline}
 \|u\|_{\tilde{L}_t^{\infty} \dot{B}_{2,1}^{s-1}}+ \|q\|_{\tilde{L}_t^{\infty} \dot{B}_{2,1}^{s-1}}+ \nu\|q\|_{\tilde{L}_t^{\infty} \dot{B}_{2,1}^{s}}+ \nu_0\|u\|_{\tilde{L}_t^1 \dot{B}_{2,1}^{s+1}}+ \nu\|q\|_{\tilde{L}_t^1 \dot{B}_{\aa}^{s+1,s-1}}+ \nu^2\|q\|_{\tilde{L}_t^1 \dot{B}_{\aa}^{s+2,s}}\\
\leq C_{p,\frac{\nu^2}{4\kappa}} e^{\displaystyle{C_{p,\frac{\nu^2}{4\kappa}} C_{visc}\int_0^t (\|\nabla v(\tau)\|_{\dot{B}_{2,1}^\fd}+ \|v(\tau)\|_{\dot{B}_{2,1}^\fd}^2)d\tau}}\\
\times\Big(\|u_0\|_{\dot{B}_{2,1}^{s-1}}+ \|q_0\|_{\dot{B}_{2,1}^{s-1}} + \nu\|q_0\|_{\dot{B}_{2,1}^{s}} +\|F\|_{\tilde{L}_t^1 \dot{B}_{2,1}^{s-1}}+ \nu\|F\|_{\tilde{L}_t^1 \dot{B}_{2,1}^{s}}+ \|G\|_{\tilde{L}_t^1 \dot{B}_{2,1}^{s-1}}\Big).
\label{estimapriori}
\end{multline}
\label{apriori}
where
$$
\begin{cases}
 \displaystyle{C_{p,\frac{\nu^2}{4\kappa}}=C \max(\sqrt{p}, \frac{1}{\sqrt{p}}) \max(\frac{4\kappa}{\nu^2}, (\frac{\nu^2}{4\kappa})^2),}\\
\displaystyle{C_{visc}=\frac{1+|\lambda+\mu|+\mu+\nu}{\nu_0}+\max(1, \frac{1}{\nu^3}).}
\end{cases}
$$
\end{thm}
\begin{rem}
 \sl{The coefficient $C_{visc}$ satisfies:
$$C_{visc}=
\begin{cases}
 \frac{1+2\nu}{\mu}+\max(1, \frac{1}{\nu^3}) & \mbox{If } \lambda+\mu>0,\\
\frac{1+2\mu}{\nu}+\max(1, \frac{1}{\nu^3}) & \mbox{If } \lambda+\mu\leq0
\end{cases}
$$
and when both viscosities are small, we simply have $C_{visc}\leq \max(1,\frac{1}{\nu_0^3})$.
}
\end{rem}

\subsection{Linear estimates}

As in \cite{CD} and \cite{CHVP} the first step to prove theorem \ref{apriori} is to obtain estimates for the following linearized system:
$$
\begin{cases}
\begin{aligned}
&\d_t q+ \div u= F,\\
&\d_t u-\cA u+ p\nabla q-\kappa \aa^2\nabla(\phia*q-q)= G.\\
\end{aligned}
\end{cases}
\leqno{(OP_\aa)}
$$
With
$$\cA u= \mu \Delta u+ (\lambda+\mu)\nabla \div u.$$
In this article, we will use the following frequency-localized estimate:
\begin{prop}
 \sl{Let $\aa>0$, $s\in \R$, $I=[0,T[$ or $[0, +\infty[$. Assume that $(q,u)$ is a solution of System $(O_\ee)$ defined on $I$. There exists $\aa_0>0$, a constant $C>0$ depending on $d$, $s$, $c_0$ and $C_0$ such that if $\aa\geq \aa_0$, for all $t\in I$ (we recall that $\nu_0$ and $C_{p,\frac{\nu^2}{4\kappa}}$ are defined in the previous theorem), and for all $j\in \Z$,
\begin{multline}
 \|\ddj u\|_{L_t^\infty L^2}+\nu_0 2^{2j} \|\ddj u\|_{L_t^1 L^2}+(1+\nu 2^j)\left(\|\ddj q\|_{L_t^\infty L^2} +\nu\min(\aa^2, 2^{2j})\|\ddj q\|_{L_t^1 L^2}\right)\\
\leq C_{p,\frac{\nu^2}{4\kappa}} \left( (1+\nu 2^j)\|\ddj q_0\|_{L^2} +\|\ddj u_0\|_{L^2} + (1+\nu 2^j)\|\ddj F\|_{L_t^1 L^2} +\|\ddj G\|_{L_t^1 L^2}\right)
\end{multline}
}
\label{estimlinloc}
\end{prop}
\begin{rem}
\sl{Let us precise that these linear estimates are valid for all dimension, the limitation $d\leq 4$ only appears in the advected case.}
\end{rem}

\subsubsection{Eigenvalues and eigenvectors}

In this article, since the methods are very close to \cite{CH} and \cite{CHVP} we will only point out what is different and refer to these articles for details. As in \cite{Dinv} or \cite{CD} we first introduce the Helmholtz decomposition of $u$. Defining the pseudo-differential operator $\Lambda$ by $\Lambda f=\mathcal{F}^{-1}(|.|\hat{f}(.))$, we set:
\begin{equation}
\begin{cases}
 v=\Lambda^{-1} \div u,\\
w=\Lambda^{-1} \curl u\\
\end{cases} 
\label{eqhelmoltz}
\end{equation}
then $u=-\Lambda^{-1}\nabla v+\Lambda^{-1}\div w$ and the system turns into:
$$
\begin{cases}
\begin{aligned}
&\d_t q+\Lambda v=F,\\
&\d_t v-\nu \Delta v -p \Lambda q +\kappa \aa^2 \Lambda(\phi_\aa*q-q)=\Lambda^{-1} \div G,\\
&\d_t w-\mu \Delta w=\Lambda^{-1} \curl G.
\end{aligned}
\end{cases}
\leqno{(L_\aa')}
$$
The last equation is a decoupled heat equation, easily estimated in Besov spaces (see \cite{Dbook} chapter 2). Moreover, as the external forces appear through homogeneous pseudo-differential operators of degree zero, we can compute the estimates in the case $F=G=0$ and deduce the general case from the Duhamel formula. So we can focus on the first two lines and compute the eigenvalues and eigenvectors of the matrix associated to the Fourier transform of the system:
$$
\d_t\left(\begin{array}{c}\hat{q}\\ \hat{v}\end{array}\right)
=A(\xi)\left(\begin{array}{c}\hat{q}\\\hat{v}\end{array}\right)\quad\hbox{with}\quad
A(\xi):= \left(\begin{array}{cc}0&-|\xi|\\|\xi|(p+\kappa \frac{|\xi|^2}{\frac{|\xi|^2}{\aa^2}+1})&-\nu|\xi|^2\end{array}\right).
$$
The discriminant of the characteristic polynomial of $A(\xi)$ is:
$$
\Delta(\xi)=|\xi|^2\left(\nu^2|\xi|^2-4(p+\kappa \frac{|\xi|^2}{\frac{|\xi|^2}{\aa^2}+1})\right),
$$
and thanks to the variations of function
$$
f_\aa:x\mapsto \nu^2 x-4(p+\kappa \frac{x}{\frac{x}{\aa^2}+1})=\nu^2 x-4(p+\kappa \aa^2) + \frac{4\kappa \aa^2}{\frac{x}{\aa^2}+1},
$$
we obtain the existence of a unique threshold $x_\aa>0$ such that
$$
\Delta(\xi)\begin{cases}
<0 \mbox{ if } |\xi|^2<x_\aa,\\
>0 \mbox{ if } |\xi|^2>x_\aa.
\end{cases}
$$
We emphasize that this function has the same variations as in the case of \cite{CHVP}: when $\frac{\nu^2}{4K}\geq 1$, $f_\aa$ is an increasing function on $\R_+$, and when $\frac{\nu^2}{4K}< 1$, $f_\aa$ is decreasing in $[0, \aa^2(\frac{2\sqrt{K}}{\nu}-1)]$ and then increasing.
\begin{prop}
 \sl{Under the same assumptions, we have:
$$
x_\aa \underset{\aa\rightarrow \infty}{\sim}
\begin{cases}
\vspace{0.2cm}
\displaystyle{\frac{4p}{\nu^2-4\kappa}} & \mbox{if }\frac{\nu^2}{4\kappa}>1,\\
\vspace{0.2cm}
\displaystyle{\aa\sqrt{\frac{p}{\kappa}}} & \mbox{if }\frac{\nu^2}{4\kappa}=1,\\
\displaystyle{(\frac{4\kappa}{\nu^2}-1)\aa^2} & \mbox{if }\frac{\nu^2}{4\kappa}<1.\\
\end{cases}
$$
}
\end{prop}
\textbf{Proof: } It is simpler than in \cite{CHVP} because here we have explicit expressions for the threshold: if we put $A=\nu^2-\kappa-4p/\aa^2$, then
$$
x_\aa=\frac{\aa^2}{2\nu^2} \left(-A+\sqrt{16p\frac{\nu^2}{\aa^2}+A}\right).
$$
Next introducing the following function, we obtain the expressions of $\hat{q}$ and $\hat{v}$ exactly as in \cite{CHVP}:
\begin{equation}
 g_\aa(x)=\frac{f_\aa(x)}{\nu^2 x}=1-\frac{4}{\nu^2 x}(p+\kappa \frac{x}{\frac{x}{\aa^2}+1})
\label{ga}
\end{equation}
\textbf{-For the low frequencies ($\Delta<0$)}, when $|\xi|<\sqrt{x_\aa}$, we have:
$$
\begin{cases}
\vspace{0.2cm}
 \hat{q}(\xi)=\frac{1}{2}\left((1+\frac{i}{S(\xi)})e^{t\lambda_+} +(1-\frac{i}{S(\xi)})e^{t\lambda_-}\right)\hat{q_0}(\xi)-i\frac{e^{t\lambda_+}-e^{t\lambda_-}}{\nu|\xi|S(\xi)}\hat{v_0}(\xi),\\
\hat{v}(\xi)=i\left(p+\kappa \frac{|\xi|^2}{\frac{|\xi|^2}{\aa^2}+1}\right)\frac{e^{t\lambda_+}-e^{t\lambda_-}}{\nu|\xi|S(\xi)} \hat{q_0}(\xi)+\frac{1}{2}\left((1-\frac{i}{S(\xi)})e^{t\lambda_+} +(1+\frac{i}{S(\xi)})e^{t\lambda_-}\right)\hat{v_0}(\xi),
\end{cases}
$$
with:
\begin{equation}
S(\xi)=\sqrt{-g_\aa|\xi|^2)}=\sqrt{\frac{4}{\nu^2 |\xi|^2}(p+\kappa \frac{|\xi|^2}{\frac{|\xi|^2}{\aa^2}+1})-1}
\label{formulelow}
\end{equation}
and
$$
\lambda_{\pm}=-\frac{\nu|\xi|^2}{2}(1\pm i S(\xi)).
$$

\textbf{-For the high frequencies ($\Delta>0$)}, when $|\xi|>\sqrt{x_\aa}$, we have:
$$
\begin{cases}
\vspace{0.2cm}
 \hat{q}(\xi)=\frac{1}{2}\left((1-\frac{1}{R(\xi)})e^{t\lambda_+} +(1+\frac{1}{R(\xi)})e^{t\lambda_-}\right)\hat{q_0}(\xi)+\frac{e^{t\lambda_+}-e^{t\lambda_-}}{\nu|\xi|R(\xi)}\hat{v_0}(\xi),\\
\hat{v}(\xi)=-\left(p+\kappa \frac{|\xi|^2}{\frac{|\xi|^2}{\aa^2}+1}\right)\frac{e^{t\lambda_+}-e^{t\lambda_-}}{\nu|\xi|R(\xi)} \hat{q_0}(\xi)+\frac{1}{2}\left((1+\frac{1}{R(\xi)})e^{t\lambda_+} +(1-\frac{1}{R(\xi)})e^{t\lambda_-}\right)\hat{v_0}(\xi),
\end{cases}
$$
with:
\begin{equation}
R(\xi)=\sqrt{g_\aa(|\xi|^2)}=\sqrt{1-\frac{4}{\nu^2 |\xi|^2}(p+\kappa \frac{|\xi|^2}{\frac{|\xi|^2}{\aa^2}+1})}
\label{formulehigh}
\end{equation}
and
$$
\lambda_{\pm}=-\frac{\nu|\xi|^2}{2}(1\pm R(\xi)).
$$
\begin{rem}
 \sl{As in \cite{CHVP} it is crucial for the time integration to observe that
$$
p+\kappa \frac{|\xi|^2}{\frac{|\xi|^2}{\aa^2}+1}=\frac{\nu^2|\xi|^2}{4}(1-R(\xi))(1+R(\xi)).
$$
}
\label{remvelocity}
\end{rem}

\subsubsection{Thresholds}

As in \cite{CHVP} we can find another threshold frequency $y_\aa> x_\aa$ of size $\aa^2$ (in each case for $\frac{\nu^2}{4\kappa}$) that will enable us to push the parabolic regularization until frequencies of size $\aa$. In the present paper we will have explicit expressions. Let us first remark that rewriting function $g_\aa$ into the following form immediately implies that this is an increasing function from $[0,\infty[$ to $]-\infty, 1[$:
\begin{equation}
 g_\aa(x)=1-\frac{4p}{\nu^2} \frac{1}{x} -\frac{\aa^2}{M} \frac{1}{x+\aa^2} \quad \mbox{with } M\overset{def}{=} \frac{\nu^2}{4\kappa}
\label{gaM}
\end{equation}
If $\beta\in[0,1[$ we easily compute that there is a unique positive solution of the equation $g_\aa(x)=\beta$ given by:
$$
x_{\aa,\beta}=\frac{1}{2}\left(-A+\sqrt{\frac{16p}{(1-\beta)\nu^2}\aa^2+A^2}\right)\quad \mbox{where } A=\frac{1}{M}\aa^2 (M-\frac{1}{1-\beta})-\frac{1}{1-\beta} \frac{4p}{\nu^2}.
$$
so that we immediately have the following asymptotics when $\aa$ is large:
$$
x_{\aa,\beta}\sim
\begin{cases}
\vspace{0.2cm}
\displaystyle{\left(\frac{1}{M(1-\beta)}-1\right) \aa^2} & \mbox{if } M<\frac{1}{1-\beta},\\
\vspace{0.2cm}
\displaystyle{\frac{2}{\nu}\sqrt{\frac{p}{1-\beta}} \aa} & \mbox{if } M=\frac{1}{1-\beta},\\
\displaystyle{\frac{4p}{\nu^2} \frac{M}{(1-\beta)M-1}} & \mbox{if } M>\frac{1}{1-\beta}.\\
\end{cases}
$$
\begin{rem}
\sl{The previous proposition is obviously a particular case of this result.}
\end{rem}
We are now able to define the second threshold $y_\aa$:
\begin{equation}
g_\aa(y_\aa)=
\begin{cases}
\vspace{0.2cm}
\displaystyle{\frac{1}{2} \mbox{ if } M\leq 1,}\\
\displaystyle{1-\frac{1}{2M}\geq \frac{1}{2} \mbox{ if } M\geq 1,}
      \end{cases}
\mbox{where }M=\frac{\nu^2}{4\kappa}.
 \end{equation}
Using the previous result for $\beta=\frac{1}{2}$ or $1-\frac{1}{2M}$ according to the case for $M$, we obtain
\begin{prop}
\sl{
With the same notations we have that:
\begin{itemize}
\item If $M\geq 1$, $y_\aa\underset{\aa \rightarrow \infty}{\sim} \aa^2$ and for all $\aa$, 
$$\aa^2 \leq y_\aa \leq 2 \aa^2.$$
\item If $M\leq 1$, $y_\aa\underset{\aa \rightarrow \infty}{\sim} (\frac{2}{M}-1)\aa^2$ and for all $\aa$,
$$\aa^2 \leq (\frac{2}{M}-1)\aa^2 \leq y_\aa \leq (\frac{2}{M}-\frac{1}{2}) \aa^2.$$
\end{itemize}
}
\end{prop}

\subsubsection{Pointwise estimates}

Now that we have defined the frequency thresholds $x_\aa$ and $y_\aa$ we have the following estimates. Up to the values of $m$ and the second exponential from the density in the second case, they are the same as in \cite{CHVP} to where we refer for details or proofs:

\begin{prop}
 \sl{Under the previous notations, there exists a constant $C$, such that for all $j\in \Z$ and all $\xi\in 2^j \cC$ where $\cC$ is the annulus $\{\xi\in \R^d, c_0=\frac{3}{4}\leq |\xi|\leq C_0=\frac{8}{3}\}$, we have the following estimates (we denote by $f_j=\ddj f$ and we refer to the appendix for details on the Littlewood-Paley theory):
\begin{itemize}
 \item If $|\xi|<\sqrt{x_\aa}$:
$$
\begin{cases}
\vspace{0.2cm}
 (1+\nu 2^j)|\hat{q_j}(\xi)|\leq C e^{-\frac{\nu t c_0^2 2^{2j}}{4}} \left((1+\nu 2^j)|\hat{q_{0,j}}(\xi)|+(1+\frac{1}{\sqrt{p}})|\hat{v_{0,j}}(\xi)|\right),\\
|\hat{v_j}(\xi)|\leq C e^{-\frac{\nu t c_0^2 2^{2j}}{4}} \left((1+\nu 2^j)(1+\sqrt{p})(1+\frac{4\kappa}{\nu^2})|\hat{q_{0,j}}(\xi)|+|\hat{v_{0,j}}(\xi)|\right).
\end{cases}
$$
\item If $\sqrt{x_\aa}<|\xi|<\sqrt{y_\aa}$:
$$
\begin{cases}
\vspace{0.2cm}
 (1+\nu 2^j)|\hat{q_j}(\xi)|\leq \frac{C}{1-m} e^{-\frac{\nu t c_0^2 2^{2j}}{4}(1-m)} \left((1+\nu 2^j)|\hat{q_{0,j}}(\xi)|+(1+\frac{1}{\sqrt{p}})|\hat{v_{0,j}}(\xi)|\right),\\
|\hat{v_j}(\xi)|\leq \frac{C}{1-m} e^{-\frac{\nu t c_0^2 2^{2j}}{4}(1-m)} \left(\nu 2^j|\hat{q_{0,j}}(\xi)|+|\hat{v_{0,j}}(\xi)|\right),
\end{cases}
$$
where $m=\sqrt{g_\aa(y_\aa)}=\frac{1}{\sqrt{2}}$ if $M=\frac{\nu^2}{4\kappa}\leq 1$, $m=\sqrt{1-\frac{1}{2M}}$ if $M\geq 1$.
 \item If $|\xi|>\sqrt{y_\aa}>\sqrt{x_\aa}$:
$$
\begin{cases}
\vspace{0.2cm}
 (1+\nu 2^j)|\hat{q_j}(\xi)|\leq C \left(e^{-\frac{\nu t |\xi|^2}{2}} +e^{-\frac{\kappa}{2\nu} \aa^2 t}\right) \left((1+\nu 2^j)|\hat{q_{0,j}}(\xi)|+(1+\frac{1}{\sqrt{p}})|\hat{v_{0,j}}(\xi)|\right),\\
|\hat{v_j}(\xi)|\leq C \left(e^{-\frac{\nu t c_0^2 2^{2j}}{4}} +\big(1-\sqrt{g_\ee(c_0^2 2^{2j})}\big)e^{-\frac{\nu t c_0^2 2^{2j}}{2}\big(1-\sqrt{g_\ee(C_0^2 2^{2j})}\big)} \right)\left(\nu 2^j |\hat{q_{0,j}}(\xi)|+|\hat{v_{0,j}}(\xi)|\right).
\end{cases}
$$
\end{itemize}
\label{Estimxi}
}
\end{prop}

\subsubsection{Time estimates}

As in \cite{CHVP}, due to the choice $c_0=3/4$ and $C_0=8/3$ (see the appendix), we can observe that there exist at most two indices $\underline{j}_\aa =\overline{j}_\aa-1$ or $\underline{j}_\aa =\overline{j}_\aa$  such that $\sqrt{y_\aa} \in 2^j[c_0, C_0]$ for $j\in\{\underline{j}_\aa,\overline{j}_\aa\}$.

We refer to \cite{CHVP} for the proof of the following proposition that implies Proposition \ref{estimlinloc}.

\begin{prop}
 \sl{Under the same assumptions as in Proposition \ref{estimlinloc}, there exists a constant $C$ such that for all $j\in \Z$ (denoting $M=\frac{\nu^2}{4\kappa}$):
\begin{itemize}
 \item For all $j\leq \overline{j}_\aa$,
\begin{multline}
\|v_j\|_{L_t^\infty L^2}+ \nu 2^{2j}\|v_j\|_{L_t^1 L^2} +(1+\nu2^j)\left(\|q_j\|_{L_t^\infty L^2}+ \nu 2^{2j}\|q_j\|_{L_t^1 L^2}\right) \leq \\
C \max(\frac{1}{M},M^2) \left((1+\nu 2^j)(1+\sqrt{p})\|q_{0,j}\|_{L^2} +(1+\frac{1}{\sqrt{p}}) \|v_{0,j}\|_{L^2}\right),\\
\label{estimlowint}
\end{multline}
\item For all $j>\overline{j}_\aa$,
\begin{multline}
\|v_j\|_{L_t^\infty L^2}+ \nu 2^{2j}\|v_j\|_{L_t^1 L^2} +(1+\nu2^j)\left(\|q_j\|_{L_t^\infty L^2}+ \frac{\nu}{\ee^2}\|q_j\|_{L_t^1 L^2}\right) \leq \\
C \max(1,M) \left((1+\nu 2^j)\|q_{0,j}\|_{L^2} +(1+\frac{1}{\sqrt{p}}) \|v_{0,j}\|_{L^2}\right).\\
\label{estimhighint}
\end{multline}
\end{itemize}
}
\label{estimlinseuil}
\end{prop}

\subsection{Advected linear estimates}

The difficulties and methods exposed here are the same as in \cite{CHVP}, so we will roughly explain them and focus on what is new.

In order to prove Theorem \ref{apriori}, a natural idea is to use Proposition \ref{estimlinloc} and put the advection terms as external forces. Unfortunately, there are some obstacles: the main problem is that in $v\cdot \nabla q$, the term $\dot{T}_v \nabla q$ can be estimated in $\dot{B}_{2,1}^{s-1}$ but not in $\dot{B}_{2,1}^s$ because it is not enough regular in high frequencies.

A direct use of the linear estimates will be useful only for the low frequencies ($j\leq 0$), and in the high frequency regime ($j>0$), we will perform a Langrangian change of variable (as in \cite{TH1}, \cite{TH2}, \cite{Dlagrangien}, \cite{CD}, \cite{CHVP}) in order to get rid of $v\cdot \nabla q$. We then aim to use on the new system our linear estimates but we have to be careful with the external force terms introduced by the change of variable. Most of the work in \cite{CHVP} was to provide estimates on the commutator of the non-local operator from the capillarity term and the Lagrangian change of variable.

For all $j\in \Z$ and $t\in I$, we introduce:
\begin{equation}
 U_j(t)=\|\ddj u\|_{L_t^\infty L^2}+\nu_0 2^{2j} \|\ddj u\|_{L_t^1 L^2}+(1+\nu 2^j)\left(\|\ddj q\|_{L_t^\infty L^2} +\nu\min(\frac{1}{\ee^2}, 2^{2j})\|\ddj q\|_{L_t^1 L^2}\right)
\label{Uj}
\end{equation}
and
\begin{multline}
 U(t)= \|u\|_{\tilde{L}_t^{\infty} \dot{B}_{2,1}^{s-1}}+ \|q\|_{\tilde{L}_t^{\infty} \dot{B}_{2,1}^{s-1}}+ \nu\|q\|_{\tilde{L}_t^{\infty} \dot{B}_{2,1}^{s}}+ \nu_0\|u\|_{\tilde{L}_t^1 \dot{B}_{2,1}^{s+1}}+ \nu\|q\|_{\tilde{L}_t^1 \dot{B}_{\ee}^{s+1,s-1}}+ \nu^2\|q\|_{\tilde{L}_t^1 \dot{B}_{\ee}^{s+2,s}}.
\end{multline}

\subsubsection{Low frequencies}
For the low frequencies, we obtain (see \cite{CHVP} section $3.1$ for details) that there exists a nonnegative summable sequence whose sum is $1$, denoted by $(c_j(t))_{j\in\Z}$ such that for all $j\leq 0$, $K>0$ (to be chosen later), and if $\aa\geq N_1/\log 2$ (we refer to \cite{CHVP} section $3.1$ for this, and to the appendix for $N_1$ which is a constant related to $c_0$ and $C_0$ in the Littlewood-Paley decomposition, such that if $|j-l|>N_1$ then $\ddj \circ \ddl =0$, in our choice of $c_0$ and $C_0$, $N_1=1$).
\begin{multline}
U_j(t) \leq C_{p,\frac{\nu^2}{4\kappa}} \Bigg[U_j(0) + (1+\nu 2^j)\|\ddj F\|_{L_t^1 L^2}+ \|\ddj G\|_{L_t^1 L^2}\\
+\frac{1}{2K} 2^{-j(s-1)}\int_0^t c_j(\tau) \left(\nu_0 \|u\|_{\dot{B}_{2,1}^{s+1}}+ \nu \|q\|_{\dot{B}_{\aa}^{s+1,s-1}}+ \nu^2\|q\|_{\dot{B}_{\aa}^{s+2,s}} \right) d\tau\\
+C^2\frac{K}{2} 2^{-j(s-1)}\int_0^t c_j(\tau) \left( \big(\max(1,\frac{1}{\nu^3}) +\frac{1}{\nu_0}\big) \|v(\tau)\|_{\dot{B}_{2,1}^\fd}^2
+C\|v(\tau)\|_{\dot{B}_{2,1}^{\fd+1}}\right) U(\tau) d\tau\Bigg].
\label{energieBF2}
\end{multline}
\begin{rem}
\sl{As pointed out in \cite{CHVP} (remark $29$) using the linear estimates for the low frequency case, allows us to get rid of another difficulty introduced by the change of variables: in low frequencies some of the additional external force terms have too much regularity to be absorbed by the left-hand side, and not enough regularity to be controlled with a view to apply the Gronwall lemma. The only way to control them would be to use interpolation arguments that would introduce linear time dependant coefficients, and prevent us to get global in time results.
}
\end{rem}

\subsubsection{Lagrangian change of coordinates}

As explained, in order to get rid of the advection terms involved in system $(LOP_\aa)$ the first step is to consider the following localized equations (as usual we set $f_j=\ddj f$...): 
$$
\begin{cases}
\begin{aligned}
&\d_t q_j+ \dot{S}_{j-1}v.\nabla q_j+ \div u_j= f_j,\\
&\d_t u_j+ \dot{S}_{j-1}v.\nabla u_j -\cA u_j+ p\nabla q_j-\kappa \aa^2 \nabla(\phia*q_j-q_j)= g_j,\\
\end{aligned}
\end{cases}
$$
where the external force terms are defined by:
$$
f_j= F_j +\left(\dot{S}_{j-1}v.\nabla q_j-\ddj (v.\nabla q)\right)\mbox{ and }g_j= G_j +\left(\dot{S}_{j-1}v.\nabla u_j-\ddj (v.\nabla u)\right).
$$
Both of these terms can be estimated thanks to the following commutator estimate from \cite{Dlagrangien} (we refer to lemma $B.1$ from appendix $B$):
\begin{lem} (\cite{Dlagrangien})
\sl{
There exists a sequence $(c_j)_{j\in \Z} \in l^1(\Z)$ such that $\|c\|_{l^1(\Z)}=1$ and a constant $C=C(d,\sigma)$ such that for all $j\in \Z$,
$$
\|\dot{S}_{j-1}v.\nabla h_j-\ddj (v.\nabla h)\|_{L^2} \leq C c_j 2^{-j\sigma}\|\nabla v\|_{\dot{B}_{2,\infty}^{\fd}\cap L^\infty} \|h\|_{\dot{B}_{2,1}^\sigma}
$$
}
\label{lemmeB}
\end{lem}
In order to perform the change of variable we define $\psi_{j,t}$ as the flow associated to $\dot{S}_{j-1}v$:
\begin{equation}
 \begin{cases}
 \partial_t \psi_{j,t}(x)=\dot{S}_{j-1}v(t,\psi_{j,t}(x))\\
\psi_{j,0}(x)=x.
\end{cases}
\label{flotdef}
\end{equation}
we can also write:
$$
\psi_{j,t}(x)=x+\int_0^t \dot{S}_{j-1}v(\tau,\psi_{j,\tau}(x)) d\tau.
$$
Thanks to propositions \ref{p:flow} and \ref{detjacobien} from the appendix (we refer to \cite{Dlagrangien} or \cite{CD}), there exists a constant $C$ such that:
\begin{equation}
 \begin{array}{lll}
\|g\circ\psi_{j,t}\|_{L^p}&\leq& e^{CV}\|g\|_{L^p}\quad\hbox{for all function }\ g\ \hbox{ in }\ 
L^p,\\[1ex]
\|D\psi_{j,t}^\pm\|_{L^\infty}&\leq&e^{CV},\\[1ex]
\|D\psi_{j,t}^\pm-I_d\|_{L^\infty}&\leq& e^{CV}-1,\\[1ex]
\|D^k\psi_{j,t}^\pm\|_{L^\infty}&\leq& C2^{(k-1)j}\Bigl(e^{CV}-1\Bigr)\quad\hbox{for }\ 
k\geq 2,
\end{array}
\label{estimflot}
\end{equation}
where
\begin{equation}
 V(t)\overset{def}{=}\int_0^t \|\nabla v(\tau)\|_{L^\infty}d\tau. 
\label{defV}
\end{equation}
As explained in \cite{CHVP}, considering some of the additionnal external force terms at the point $\psi_{j,t}^{-1}(x)$ instead of $x$ will help getting uniform estimates with respect to $\aa$, so the jacobian determinant of the change of variable will play a crucial role in our computations (contrary to the case of lemma $2.6$ from \cite{Dbook} where it produces a term that we are not able to sum here).
\begin{equation}
\begin{cases}
det(D\psi_{j,t}(x))=e^{\int_{0}^t (\div \dot{S}_{j-1}v)(\tau,\psi_{j,\tau(x)})d\tau},\\
det(D\psi_{j,t}^{-1}(x))=e^{-\int_{0}^t (\div \dot{S}_{j-1}v)(\tau, X_j(\tau,t,x))d\tau} =e^{-\int_{0}^t (\div \dot{S}_{j-1}v)(\tau,\psi_{j,\tau}\circ \psi_{j,t}^{-1}(x))d\tau},
\end{cases}
\label{detjacobienj}
\end{equation}
where $X_j(\tau,t,x))$ denotes the two parameter flow associated to $\dot{S}_{j-1}v$ (we refer to \eqref{flot2param} in the appendix).

Let us now perform the lagrangian change of variable, for a function $h$, we define $\tilde h =h \circ \psi_{j,t} =h(t,\psi_{j,t})$. Then we have $\partial_t \tilde q_j(t,x)=(\partial_t q_j+ \dot{S}_{j-1} v.\nabla q_j)(t, \psi_{j,t}(x))$, which provides the following system:
\begin{equation}
 \begin{cases}
\begin{aligned}
&\d_t \tilde q_j+ \div \tilde u_j= \tilde f_j +R_j^1,\\
&\d_t \tilde u_j -\cA \tilde u_j+ p\nabla \tilde q_j-\kappa \aa^2 \nabla(\phia *\tilde q_j-\tilde q_j)= \tilde g_j +R_j^2 +R_j^3 +\kappa R_j,\\
\end{aligned}
\end{cases}
\label{systchanged}
\end{equation}
where the remainder terms $R_q^1,$ $R_q^2$ and $R_q^3$ are exactly the same as in \cite{CD} and \cite{CHVP} (with the same convention: if $f:\R^d\rightarrow \R^m$ is a differentiable function then $Df$ denotes the Jacobian matrix of $f,$ and $\nabla f$ is the transposed matrix of $Df.$):
$$\displaylines{
R_j^1(t,x):={\rm Tr}\bigl(\nabla \tilde u_j(t,x)\cdot (I_d-\nabla\psi_{j,t}^{-1}(\psi_{j,t}(x)))\bigr),\cr
R_j^2(t,x):= \nabla \tilde q_j(t,x)\cdot (\nabla\psi_{j,t}^{-1}(\psi_{j,t}(x))-I_d)}
$$
and  $R_j^3:= \mu R_j^4 + (\lambda+\mu) R_j^5$ with
$$
\displaylines{
\quad R_j^{4,i}(t,x):={\rm Tr}\Bigl( (\nabla\psi_{j,t}^{-1}(\psi_{j,t}(x))-I_d)\cdot\nabla  D \tilde u_j^i(t,x) \cdot D\psi_{j,t}^{-1}(\psi_{j,t}(x))\hfill \cr\hfill
+\nabla D\tilde u_j^i(t,x)\cdot (D\psi_{j,t}^{-1}(\psi_{j,t}(x))-I_d)\Bigr) +\nabla \tilde u_j^i(t,x)\cdot 
\Delta \psi_{j,t}^{-1}(\psi_{j,t}(x)) \cr
\quad R_j^{5,k}(t,x):={\rm Tr}\Bigl( D \tilde u_j(t,x) \cdot \partial_k D\psi_{j,t}^{-1}(\psi_{j,t}(x))\Bigr)\hfill\cr\hfill
+ \Sum_{a,b,c, b\neq a, c\neq k} \partial_{bc}^2 \tilde u_j^i(t,x)\cdot \partial_k \psi_{j,t}^{-1,c}(\psi_{j,t}(x)) \cdot \partial_a \psi_{j,t}^{-1,b}(\psi_{j,t}(x))\hfill \cr\hfill
+ \Sum_{i=1}^d \partial_{ki}^2 \tilde u_j^i(t,x)\cdot \Bigl(\partial_k \psi_{j,t}^{-1,k}(t,\psi_{j,t}(x))-I_d) \cdot \partial_i \psi_{j,t}^{-1,i}(t,\psi_{j,t}(x))+ (\partial_i \psi_{j,t}^{-1,i}(\psi_{j,t}(x))-I_d) \Bigr).}
$$
As in \cite{CHVP}, the only difference with \cite{CD} is the following additionnal remainder term:
\begin{equation}
R_j=\aa^2 (\phia*\nabla q_j-\nabla q_j) \circ \psi_{j,t} -\aa^2 (\phia*\nabla \tilde q_j-\nabla \tilde q_j).
\label{Rj}
\end{equation}
Thanks to the definition of $\phia$ and $\phi$ we can also write that:
\begin{equation}
 L_\aa(f)\overset{def}{=}\aa^2 (\phia* f-f)= \aa^2 \int_{\R^d} \phi(z)\left(f(x-\frac{z}{\aa})-f(x)\right)dz.
\label{phiaa}
\end{equation}
As in \cite{CHVP} most of the work consists in obtaining bounds in $L_t^1 L^2$ that are uniform with respect to $\aa$, and go to zero when $t$ is small and dealing with $R_j$ is the object of the rest of this section.

\subsubsection{Precisions on the capillary term}

Let us first go back to the convolution term written in (\ref{phiaa}): for a function $f$,
$$
\aa^2 \hat{(\phi_\ee* f-f)}(\xi)=-\frac{|\xi|^2}{1+\frac{|\xi|^2}{\aa^2}} \hat{f}(\xi).
$$
Similarly to \cite{CHVP} we obtain the following equivalence, giving a smooth interpretation of the hybrid norm. Here again, instead of a fixed frequency threshold there is a continuous transition from the parabolically regularized low frequencies and the damped high frequencies:
\begin{prop}
 \sl{For any suitable function $f$ and any $s\in\R$, we have:
\begin{multline}
 \|f\|_{\dot{B}_\aa^{s+2,s}} =\sum_{j\in \Z} \min(\aa^2, 2^{2j}) 2^{js} \|\ddj f\|_{L^2}\\
 \sim \sum_{j\in \Z} \frac{2^{2j}}{1+\frac{2^{2j}}{\aa^2}}  2^{js} \|\ddj f\|_{L^2} \sim \|\aa^2 (\phi_\ee*f-f)\|_{\dot{B}_{2,1}^s}
\label{normhybride2}
\end{multline}
}
\end{prop}
\textbf{Proof:} Thanks to the monotonicity of function $x\mapsto x/(1+x)$, we easily obtain that for all $j\in \Z$:
$$
\frac{1}{2} \min (\aa^2,2^{2j}) \leq \frac{2^{2j}}{1+\frac{2^{2j}}{\aa^2}} \leq \min (\aa^2,2^{2j}).
$$
We refer to \cite{CHVP} for details.

\begin{rem}
 \sl{In the $L^r$-setting, we can prove that for all $j\in\Z$:
$$
\|\aa^2(\phi_\ee*\ddj f-\ddj f)\|_{L^r} \leq C\min(\aa^2, 2^{2j}) \|\ddl f\|_{L^r}.
$$
}
\label{Lr}
\end{rem}

\subsubsection{Estimates on the capillary remainder $R_j$}
\label{sectionRj}
This section is devoted to giving estimates on the capillary term introduced in \eqref{Rj}. As
$$
\nabla \tilde q_j =\nabla( q_j\circ \psi_{j,t})=\nabla q_j \circ \psi_{j,t} \times D\psi_{j,t} =\nabla q_j \circ \psi_{j,t} \times (D\psi_{j,t}-I_d) +\nabla q_j \circ \psi_{j,t},
$$
we obtain the following decomposition: $R_j= I_j+II_j$ with
\begin{equation}
\begin{cases}
\vspace{0.2cm}
I_j= \aa^2(\phi_\ee*g_j-g_j) \mbox{ where } g_j= \nabla q_j \circ \psi_{j,t} \times (I_d-D\psi_{j,t}),\\
II_j= \aa^2(\phi_\ee*\nabla q_j-\nabla q_j)\circ \psi_{j,t} -\aa^2\Big(\phi_\ee*(\nabla q_j \circ \psi_{j,t}) -\nabla q_j \circ \psi_{j,t}\Big).
\end{cases}
\label{IjIIj}
\end{equation}
So that as in \cite{CHVP}, we need to estimate (locally in frequency) the commutator between the Lagrangian change of variable and the non-local operator $L_\aa$ defined in \eqref{phiaa}. For a function $f$, and all $j\in\Z$ we set $f_j=\ddj f$ and:
\begin{equation}
 II_j'= II_j'(f)= \aa^2(\phia *f_j-f_j)\circ \psi_{j,t} -\aa^2\Big(\phia *(f_j \circ \psi_{j,t}) -f_j \circ \psi_{j,t}\Big).
\end{equation}

\begin{thm}
 \sl{Let $\sigma \in \R$. There exists a constant $C=C_{\sigma, d}$ such that for all $f\in\dot{B}_\aa^{\sigma+2, \sigma}$, there exists a summable positive sequence $(c_j(f))_{j\in \Z}$ whose sum is $1$ such that for all $t$ so small that
\begin{equation}
e^{2CV}-1 \leq \frac{1}{2}.
 \label{Cond1}
\end{equation}
and for all $j\in \Z$,
$$
\|II_j'(f)\|_{L^2}\leq C e^{CV}(V+e^{2CV}-1) c_j(f) 2^{-j\sigma}\|\aa^2(\phia*f-f)\|_{\dot{B}_{2,1}^\sigma},
$$
where $V(t)=\int_0^t \|\nabla v(\tau)\|_{L^\infty}d\tau$. 
}
\label{II'}
\end{thm}
\begin{rem}
\sl{As a by-product we obtain that under the previous assumptions, if $t$ is small enough,
$$
\aa^2 \|(\phi_\ee*\ddj f)\circ \psi_{j,t} -\phi_\ee*(\ddj f\circ \psi_{j,t})\|_{L^2} \leq C e^{CV}(V+e^{2CV}-1) c_j(f) 2^{-j\sigma}\|\aa^2(\phi_\ee*f-f)\|_{\dot{B}_{2,1}^\sigma}.
$$
Remark that neither of the left-hand side terms are spectrally localized.} 
\end{rem}
\textbf{Proof:} we refer to \cite{CHVP} for details and here we will only focus on what changes. The first step is to obtain pointwise estimates and then $L^2$ estimates: as in the works of T. Hmidi, S. Keraani, H. Abidi and M. Zerguine (\cite{TH1}, \cite{TH2}, \cite{TH3} and \cite{TH4}) we wish to retrieve the desired Besov norm thanks to an equivalent expression of $II_j'$ as an integral formulation involving finite differences of $f$ of order $1$ (that is expressions of the type $\tau_{-y} f - f$ where $\tau_{-y} f(x)= f(x+y)$) or order $2$. Like in \cite{CHVP} we need to directly consider $II_j'(\psi_{j,t}^{-1}(x))$ instead of simply $II_j'(x)$ and we obtain:
\begin{multline}
 II_j'(\psi_{j,t}^{-1}(x))= \aa^2 \int_{\R^d} \phi(z)\left(f_j(x-\frac{z}{\aa})-f_j(x)\right)\times\\
\left[1-\frac{\phi\Big(\aa \big(\psi_{j,t}^{-1}(x)-\psi_{j,t}^{-1}(x-\frac{z}{\aa})\big)\Big)}{\phi(z)} e^{\displaystyle{-\int_{0}^t (\div \dot{S}_{j-1}v)(\tau, X_j(\tau,t,x-\frac{z}{\aa}))d\tau}} \right]dz.
\end{multline}
\begin{rem}
\sl{We emphasize here that the previous quotient is well defined near zero and we refer to section \eqref{interactionpotentielgen}.}
\end{rem}
We have to face the same problem as in \cite{CHVP}: as we want to estimate the $L^2$ norm of this quantity, that is a Besov norm with integer regularity index $s=0$, the finite difference of order $1$ will not be sufficient for our need, and we will have to introduce finite differences of order $2$ (this is a classical problem for integer indices). Indeed, using the present quantity would only involve a term in $|y|/\aa$  and when estimating in low frequencies, there would be either an extra multiplicative coefficient $\aa$ or an extra derivative term $2^{-j}$ (that would prevent any convergence when $-j$ is large). To be able to do a correct estimate we need at least $|y|^2/\aa^2$.

For this, we simply write $II_j'=\frac{1}{2}(II_j'+II_j')$, and perform the change of variable $z=-y$ in the second integral. If we set:
\begin{equation}
 \begin{cases}
\vspace{0.2cm}
R_-= \displaystyle{\frac{\phi\Big(\aa \big(\psi_{j,t}^{-1}(x)-\psi_{j,t}^{-1}(x-\frac{z}{\aa})\big)\Big)}{\phi(z)},}\\
\vspace{0.2cm}
R_+= \displaystyle{\frac{\phi\Big(\aa \big(\psi_{j,t}^{-1}(x)-\psi_{j,t}^{-1}(x+\frac{z}{\aa})\big)\Big)}{\phi(z)},}\\
\vspace{0.2cm}
B_-=\displaystyle{\int_{0}^t (\div \dot{S}_{j-1}v)(\tau, X_j(\tau,t,x-\frac{z}{\aa}))d\tau,}\\
B_+=\displaystyle{\int_{0}^t (\div \dot{S}_{j-1}v)(\tau, X_j(\tau,t,x+\frac{z}{\aa}))d\tau,}
 \end{cases}
\label{RRBB}
\end{equation}
then
\begin{multline}
 II_j'(\psi_{j,t}^{-1}(x))=\frac{\aa^2}{2} \int_{\R^d} \phi(z)\left(f_j(x-\frac{z}{\aa})-f_j(x)\right)[1-R_- e^{-B_-}] dz\\
 +\frac{\aa^2}{2} \int_{\R^d} \phi(z)\left(f_j(x+\frac{z}{\aa})-f_j(x)\right)[1-R_+ e^{-B_+}] dz= III_j(x)+ IV_j(x).
\end{multline}
where
\begin{equation}
 \begin{cases}
\vspace{0.2cm}
  III_j(x)= \displaystyle{\frac{\aa^2}{2} \int_{\R^d} \phi(z)\left(f_j(x-\frac{z}{\aa})+f_j(x+\frac{z}{\aa})-2f_j(x)\right)[1-R_- e^{-B_-}] dz,}\\
IV_j(x)= \displaystyle{\frac{\aa^2}{2} \int_{\R^d} \phi(z)\left(f_j(x+\frac{z}{\aa})-f_j(x)\right)[R_- e^{-B_-}-R_+ e^{-B_+}] dz.}
 \end{cases}
\label{IIIetIV}
\end{equation}
Taking $L^2$ norms we have:
\begin{equation}
\|II_j'\circ \psi_{j,t}^{-1}\|_{L^2} \leq \|III_j\|_{L^2} +\|IV_j\|_{L^2},
\label{estimbaseIIj}
\end{equation}
and thanks to estimates on the Jacobian determinant of the flow (see \ref{estimflot}), theorem \ref{II'} is immediately implied by the following proposition:
\begin{prop}
 \sl{Under the previous assumptions, there exist a positive constant $C=C_{\sigma,d}$ and a nonnegative sequence $(c_j= c_j(f))_{j\in \Z}$ whose sum is $1$, such that if $t$ is so small that $e^{2CV(t)}-1\leq \frac{1}{2}$, we have:
$$
\|III_j\|_{L^2} +\|IV_j\|_{L^2} \leq C (V+e^{2CV}-1) e^{CV} 2^{-j\sigma} c_j \|\aa^2(\phia* f-f)\|_{\dot{B}_{2,1}^\sigma}.
$$
\label{propIII}
}
\end{prop}
To prove this result we will successively prove the following lemmas:
\begin{lem}
 \sl{There exists a constant $C$ only depending on the dimension $d$ such that for all $j\in\Z$, $f$, and all $t$ is so small that $e^{2CV(t)}-1\leq \frac{1}{2}$,
\begin{multline}
 \|III_j\|_{L^2}\leq C e^{CV} (V+e^{2CV}-1) \\
\times \aa^2 \int_{\R^d} e^{-\frac{|z|}{8}} \frac{1+|z|^{-\frac{1}{2}}}{|z|^{d-\frac{3}{2}}} \|f_j(.-\frac{z}{\aa}) +f_j(.+\frac{z}{\aa}) -2f_j(.)\|_{L^2} dz,
\end{multline}
and
\begin{multline}
 \|IV_j\|_{L^2}\leq C e^{CV} (V+e^{2CV}-1) \\
\times \aa^2 \int_{\R^d} e^{-\frac{|z|}{8}} \frac{1+|z|^{-\frac{1}{2}}}{|z|^{d-\frac{3}{2}}} \min(1,\frac{2^j|z|}{\aa}) \|f_j(.+\frac{z}{\aa})-f_j(.)\|_{L^2} dz,
\end{multline}
where $V$ is defined in \eqref{defV}.
}
\label{lemIIIa}
\end{lem}
\begin{lem}
 \sl{For all $\sigma\in \R$ there exists a constant $C_{\sigma,d}$ such that for any $f\in \dot{B}_{2,1}^\sigma$, there exists a nonnegative summable sequence $(c_j(f))_{j\in \Z}$ with $\|c_j(f)\|_{l^1(\Z)}=1$, such that
\begin{multline}
\aa^2 \int_{\R^d} e^{-\frac{|z|}{8}} \frac{1+|z|^{-\frac{1}{2}}}{|z|^{d-\frac{3}{2}}} \|f_j(.-\frac{z}{\aa}) +f_j(.+\frac{z}{\aa}) -2f_j(.)\|_{L^2} dz\\
+\aa^2 \int_{\R^d} e^{-\frac{|z|}{8}} \frac{1+|z|^{-\frac{1}{2}}}{|z|^{d-\frac{3}{2}}} \min(1,\frac{2^j|z|}{\aa}) \|f_j(.+\frac{z}{\aa})-f_j(.)\|_{L^2} dz\\
\leq C_{\sigma,d} 2^{-j\sigma} c_j(f) \|\aa^2(\phia* f-f)\|_{\dot{B}_{2,1}^\sigma}. 
\end{multline}
}
\label{lemIIIb}
\end{lem}
\textbf{Proof of lemma \ref{lemIIIa}:} A direct estimate gives:
$$
\|III_j\|_{L^2} \leq \displaystyle{\frac{\aa^2}{2} \int_{\R^d} \phi(z)\|f_j(\cdot-\frac{z}{\aa})+f_j(\cdot+\frac{z}{\aa})-2f_j(\cdot)\|_{L_x^2}\|1-R_- e^{-B_-}\|_{L_x^\infty} dz,}
$$
and
$$
\|IV_j\|_{L^2} \leq \displaystyle{\frac{\aa^2}{2} \int_{\R^d} \phi(z)\|f_j(\cdot+\frac{z}{\aa})-f_j(\cdot)\|_{L_x^2}\|R_- e^{-B_-}-R_+ e^{-B_+}\|_{L_x^\infty} dz,}
$$
So that we first have to focus on the $L^\infty$ norms:
\begin{lem}
\sl{Under the same assumptions, for all $\eta\in]0,1[$, there exists three constants $C$, $C_d$ and $C_{d,\eta}$ such that for all $j\in\Z$, all $f$, and all $t>0$ so small that $e^{2CV(t)}-1 \leq \frac{1}{2}$, we have:
$$
\|1-R_- e^{-B_-}\|_{L_x^\infty}\leq C_{d,\eta} e^{C_dV(t)} (e^{2CV}-1+V) e^{|z|\left[\eta+(1-\eta)\left(e^{2CV(t)}-1\right)\right]}  (1+|z|^{-\frac{1}{2}}).
$$
\begin{multline}
\|R_- e^{-B_-}-R_+ e^{-B_+}\|_{L_x^\infty}\\
\leq C_{d,\eta} e^{C_dV(t)} (e^{2CV}-1+V)\min(1,\frac{2^j|z|}{\aa})  e^{|z|\left[\eta+(1-\eta)\left(e^{2CV(t)}-1\right)\right]}  (1+|z|^{-\frac{1}{2}}).
\end{multline}
}
\label{estimLinf}
\end{lem}
\textbf{Proof of Lemma \ref{estimLinf}:} as in \cite{CHVP} every term is close to 1 when $t$ is small and we obtain the result by carefully estimating the differences between these terms. For this we write:
\begin{equation}
\begin{cases}
Q_1= 1-R_- e^{-B_-}=1-R_-+R_-(1-e^{-B_-}),\\
Q_2= R_- e^{-B_-}-R_+ e^{-B_+}= (R_--R_+) e^{-B_-}+R_+ (e^{-B_-}-e^{-B_+}).
\end{cases}
\label{q1q2}
\end{equation}
Thanks to the following elementary consequence of the mean-value theorem:
\begin{lem}
 \sl{For any $x,y\in \R$, $|e^x-e^y|\leq |x-y|e^{\max(x,y)}$.}
\label{lemexp}
\end{lem}
we obtain (we refer to \cite{CHVP} for details) that there exists a constant $C>0$ such that (we refer to theorem \ref{II'} for the definition of $V$):
\begin{equation}
\begin{cases}
\vspace{0.2mm}
\displaystyle{|1-e^{-B_-}|\leq C V(t) e^{CV(t)},}\\
\vspace{0.2mm}
\displaystyle{|e^{-B_-}|\leq e^{CV(t)},}\\
\displaystyle{|e^{-B_-}-e^{-B_+}| \leq 2C V(t) e^{2CV(t)} \min(1,\frac{|z|2^j}{\aa}).}
\end{cases}
\label{coeffaciles}
\end{equation}
We now turn to the estimates involving $R_\pm$. As obtained in \eqref{phiexpl}, the interaction potential can be expressed thanks to the modified Bessel function of second kind $K_{\fd-1}$: 
$$
\phi(x)=\frac{C_d}{|x|^{\fd-1}} K_{\fd-1}(|x|),
$$
so that if we denote by $Y_\pm$ the variation ratio of the flow:
\begin{equation}
Y_\pm \overset{def}{=} \frac{\big|\psi_{j,t}^{-1}(x)-\psi_{j,t}^{-1}(x\pm \frac{z}{\aa})\big|}{\frac{|z|}{\aa}}
\label{ratioY}
\end{equation}
we have the expression:
\begin{equation}
R_\pm= \left(\frac{1}{Y_\pm}\right)^{\fd-1} \frac{K_{\fd-1}(\aa \big|\psi_{j,t}^{-1}(x)-\psi_{j,t}^{-1}(x\pm \frac{z}{\aa})\big|)}{K_{\fd-1}(|z|)} =\left(\frac{1}{Y_\pm}\right)^{\fd-1} \frac{K_{\fd-1}(|z|Y_\pm)}{K_{\fd-1}(|z|)}.
\end{equation}
We bound the first factor thanks to the following estimates on the variation ratio:
\begin{lem}[\cite{CHVP} section $3.4$]
 \sl{Under the general previous assumptions, for all $x,z\in \R^d$ with $z\neq 0$, we have:
$$
e^{-C V(t)} \leq Y_\pm \leq e^{C V(t)}
$$
and
$$
\left|Y_\pm -1\right|\leq e^{2CV(t)}-1, \quad \left|\frac{1}{Y_\pm} -1\right|\leq e^{2CV(t)}-1.
$$
}
\label{estimtauxflot}
\end{lem}
Using the upper and lower bounds for $\phi$ given by proposition \ref{propK}, we can write that for a fixed $\eta\in]0,1[$ (to be precised later) there exists a constant $C_{d,\eta}$ such that:
$$
R_\pm\leq C_{d,\eta} e^{C(\fd-1)V(t)} \frac{e^{-(1-\eta)\aa \big|\psi_{j,t}^{-1}(x)-\psi_{j,t}^{-1}(x\pm \frac{z}{\aa})\big|}}{\left(\aa \big|\psi_{j,t}^{-1}(x)-\psi_{j,t}^{-1}(x\pm \frac{z}{\aa})\big|\right)^{\frac{d-1}{2}}} e^{|z|} \cdot
\begin{cases}
|z|^{\fd-1} & \mbox{ if } d>2,\\
(1+|z|^{\frac{1}{2}}) & \mbox{ if } d=2.
\end{cases}.
$$
that is, with the previous notations:
$$
R_\pm\leq C_{d,\eta} e^{C(\fd-1)V(t)} e^{|z|\left(1-(1-\eta)Y_\pm\right)}  \left(\frac{1}{Y_\pm}\right)^{\frac{d-1}{2}} \frac{1}{|z|^{\frac{d-1}{2}}}
\cdot
\begin{cases}
|z|^{\fd-1} & \mbox{ if } d>2,\\
(1+|z|^{\frac{1}{2}}) & \mbox{ if } d=2.
\end{cases}.
$$
and then
$$
R_\pm\leq C_{d,\eta} e^{C(d-\frac{3}{2})V(t)} e^{\eta |z|} e^{(1-\eta)|z|\left(1-Y_\pm\right)} (1+|z|^{-\frac{1}{2}}).
$$
Thanks again to lemma \ref{estimtauxflot}, we finally obtain:
\begin{equation}
R_\pm\leq C_{d,\eta} e^{C_dV(t)} e^{|z|\left[\eta+(1-\eta)\left(e^{2CV(t)}-1\right)\right]} (1+|z|^{-\frac{1}{2}}).
\label{estimRpm}
\end{equation}
There two terms left to estimate: $1-R_-$ and $R_--R_+$. Again, we simply write that:
$$
1-R_\pm= P_1+ P_2=\left[1-\left(\frac{1}{Y_\pm}\right)^{\fd-1}\right] +\left(\frac{1}{Y_\pm}\right)^{\fd-1} \left[1-\frac{K_{\fd-1}(\aa \big|\psi_{j,t}^{-1}(x)-\psi_{j,t}^{-1}(x\pm \frac{z}{\aa})\big|)}{K_{\fd-1}(|z|)}\right].
$$
An elementary use of the mean-value theorem to the function $h(y)=y^{\fd-1}$ gives that (we recall that in lemma \ref{estimtauxflot}, we proved $e^{-CV(t)} \leq Y_\pm \leq e^{CV(t)}$):
\begin{equation}
|h(1)-h(\frac{1}{Y_\pm})| \leq |1-\frac{1}{Y_\pm}| (\fd-1)\underset{y\in [e^{-CV(t)}, e^{CV(t)}]}{\sup}  y^{\fd-2}\leq C_d |1-\frac{1}{Y_\pm}| e^{C|\fd-2|V(t)},
\label{TVI}
\end{equation}
that is:
\begin{equation}
|P_1|\leq C_d (e^{2CV(t)}-1)e^{C_dV(t)}
\label{estimP1}
\end{equation}
Moreover, the second term is bounded by:
$$
|P_2|\leq \frac{e^{C_dV(t)}}{K_{\fd-1}(|z|)} \Big|K_{\fd-1}(|z|)-K_{\fd-1}(\aa \big|\psi_{j,t}^{-1}(x)-\psi_{j,t}^{-1}(x\pm \frac{z}{\aa})\big|)\Big|
$$
Similarly, we can write that:
\begin{multline}
\Big|K_{\fd-1}(|z|)-K_{\fd-1}(\aa \big|\psi_{j,t}^{-1}(x)-\psi_{j,t}^{-1}(x\pm \frac{z}{\aa})\big|)\Big|\\
\leq |z||Y_\pm-1| \int_0^1 \left|K_{\fd-1}'\big(|z|(1+u(|Y_\pm-1|)\big)\right| du\\
\leq (e^{2CV(t)}-1)|z|\cdot \underset{y\in \left[|z|\big(1-(e^{2CV(t)}-1)\big), |z|\big(1+(e^{2CV(t)}-1)\big)\right]}{\sup} \left|K_{\fd-1}'(y)\right|.
\label{TVIK}
\end{multline}
If $t$ is so small that $e^{2CV(t)}-1\leq \frac{1}{2}$, then thanks to the bound for $K_{\fd-1}'$ from proposition \ref{propK} we obtain that:
\begin{multline}
\Big|K_{\fd-1}(|z|)-K_{\fd-1}(\aa \big|\psi_{j,t}^{-1}(x)-\psi_{j,t}^{-1}(x\pm \frac{z}{\aa})\big|)\Big| \leq (e^{2CV(t)}-1)|z|\cdot \underset{y\in \left[\frac{1}{2}|z|, \frac{3}{2}|z|\right]}{\sup} C_{d,\eta} \frac{e^{-(1-\eta) y}}{y^{\frac{d+1}{2}}}\\
\leq C_{d,\eta} (e^{2CV(t)}-1) \frac{e^{-(1-\eta) |z|\big(1-(e^{2CV(t)}-1)\big)}}{|z|^{\frac{d-1}{2}}}.
\end{multline}
Gathering this with the lower bound for $K_{\fd-1}$ from proposition \ref{propK}, we finally bound $P_2$:
\begin{multline}
|P_2|\leq e^{C_dV(t)} C_{d,\eta} (e^{2CV(t)}-1) \frac{e^{-(1-\eta) |z|\big(1-(e^{2CV(t)}-1)\big)}}{|z|^{\frac{d-1}{2}}} \cdot e^{|z|} \cdot
\begin{cases}
|z|^{\fd-1} & \mbox{ if } d>2,\\
(1+|z|^{\frac{1}{2}}) & \mbox{ if } d=2
\end{cases}\\
\leq C_{d,\eta} e^{C_dV(t)} (e^{2CV(t)}-1) e^{\eta|z|} e^{(1-\eta) |z|(e^{2CV(t)}-1)} (1+|z|^{-\frac{1}{2}}).
\end{multline}
Finally together with \eqref{estimP1}, we finally obtain the estimate:
\begin{equation}
|1-R_\pm| \leq C_{d,\eta} e^{C_dV(t)} (e^{2CV(t)}-1) e^{|z|\left[\eta+(1-\eta)\left(e^{2CV(t)}-1\right)\right]}  (1+|z|^{-\frac{1}{2}}).
\label{estim1moinsRpm}
\end{equation}
Let us finally turn to the last term: $R_--R_+$. Using the notations, we can write:
\begin{multline}
R_--R_+=\left(\frac{1}{Y_-}\right)^{\fd-1} \frac{K_{\fd-1}(|z|Y_-)}{K_{\fd-1}(|z|)} -\left(\frac{1}{Y_+}\right)^{\fd-1} \frac{K_{\fd-1}(|z|Y_+)}{K_{\fd-1}(|z|)}\\
=\left[\left(\frac{1}{Y_-}\right)^{\fd-1} -\left(\frac{1}{Y_+}\right)^{\fd-1}\right] \frac{K_{\fd-1}(|z|Y_-)}{K_{\fd-1}(|z|)} +\left(\frac{1}{Y_+}\right)^{\fd-1} \frac{K_{\fd-1}(|z|Y_-)-K_{\fd-1}(|z|Y_+)}{K_{\fd-1}(|z|)}\\
=P_1'+P_2'.
\label{pp1pp2}
\end{multline}
As we did in \eqref{TVI}, if $k(y)=y^{-(\fd-1)}$, we can write:
$$
|k(\frac{1}{Y_-})-k(\frac{1}{Y_+})| \leq C_d |Y_--Y_+| e^{C_dV(t)},
$$
The other term in $P_1'$ has already been estimated when estimating $R_\pm$:
$$
\frac{K_{\fd-1}(|z|Y_-)}{K_{\fd-1}(|z|)} \leq C_{d,\eta} e^{C_dV(t)} e^{|z|\left[\eta+(1-\eta)\left(e^{2CV(t)}-1\right)\right]} (1+|z|^{-\frac{1}{2}}),
$$
and we get that:
\begin{equation}
|P_1'|\leq C_{d,\eta} e^{C_dV(t)} |Y_--Y_+| e^{|z|\left[\eta+(1-\eta)\left(e^{2CV(t)}-1\right)\right]} (1+|z|^{-\frac{1}{2}})
\label{estimPp1}
\end{equation}
For $P_2'$, as in \eqref{TVIK}, we can write:
\begin{multline}
\Big|K_{\fd-1}(|z|Y_-)-K_{\fd-1}(|z|Y_+)\Big| \leq |z||Y_--Y_+| \int_0^1 \left|K_{\fd-1}'\bigg(|z|\big((1-u)Y_++u Y_-\big)\bigg)\right|du\\
\leq |z||Y_--Y_+| \int_0^1 \left|K_{\fd-1}'\bigg(|z|+|z|\big((1-u)(Y_+-1)+u (Y_--1)\big)\bigg)\right|du\\
\leq |z||Y_--Y_+| \cdot \underset{y\in \left[|z|\big(1-(e^{2CV(t)}-1)\big), |z|\big(1+(e^{2CV(t)}-1)\big)\right]}{\sup} \left|K_{\fd-1}'(y)\right|\\
\leq C_{d,\eta} |Y_--Y_+| \frac{e^{-(1-\eta) |z|\big(1-(e^{2CV(t)}-1)\big)}}{|z|^{\frac{d-1}{2}}}.
\end{multline}
Using once again the lower bound for $K_{\fd-1}$ implies that:
$$
|P_2'|\leq C_{d,\eta} e^{C_dV(t)} |Y_--Y_+| e^{|z|\left[\eta+(1-\eta)\left(e^{2CV(t)}-1\right)\right]}  (1+|z|^{-\frac{1}{2}}).
$$
and together with \eqref{estimPp1}, this implies that:
\begin{equation}
|R_--R_+| \leq C_{d,\eta} e^{C_dV(t)} |Y_--Y_+| e^{|z|\left[\eta+(1-\eta)\left(e^{2CV(t)}-1\right)\right]}  (1+|z|^{-\frac{1}{2}}).
\label{estimRRbof}
\end{equation}
Thanks to lemma \ref{estimtauxflot}, we get:
\begin{equation}
|Y_--Y_+|  \leq e^{CV(t)}(e^{2CV(t)}-1).
\label{estimYYbof}
\end{equation}
Unfortunately, as explained, this will not be sufficient: this estimate is useful for high frequencies $j$, but after integration in $z$, the result is not summable when $j$ goes to $-\infty$. This is why, as in \cite{CHVP}, we need a much more precise estimate on $|Y_--Y_+|$:
\begin{multline}
|Y_--Y_+| =\frac{\aa}{|z|} \left(\big|\psi_{j,t}^{-1}(x)-\psi_{j,t}^{-1}(x -\frac{z}{\aa})\big| -\big|\psi_{j,t}^{-1}(x)-\psi_{j,t}^{-1}(x +\frac{z}{\aa})\big|\right)\\
=\frac{\aa}{|z|} \frac{\big|\psi_{j,t}^{-1}(x)-\psi_{j,t}^{-1}(x -\frac{z}{\aa})\big|^2 -\big|\psi_{j,t}^{-1}(x)-\psi_{j,t}^{-1}(x +\frac{z}{\aa})\big|^2}{\big|\psi_{j,t}^{-1}(x)-\psi_{j,t}^{-1}(x -\frac{z}{\aa})\big| +\big|\psi_{j,t}^{-1}(x)-\psi_{j,t}^{-1}(x +\frac{z}{\aa})\big|}=\frac{\aa}{|z|} \frac{N_j}{D_j}
\end{multline}
As in the proof of lemma $8$ from \cite{CHVP}, we simply use the identity $|a|^2-|b|^2=(a+b|a-b)$ and write:
\begin{multline}
N_j=\left( 2\psi_{j,t}^{-1}(x)-\psi_{j,t}^{-1}(x -\frac{z}{\aa})-\psi_{j,t}^{-1}(x +\frac{z}{\aa})\Big| \psi_{j,t}^{-1}(x -\frac{z}{\aa}) -\psi_{j,t}^{-1}(x +\frac{z}{\aa})\right)\\
\leq \big| 2\psi_{j,t}^{-1}(x)-\psi_{j,t}^{-1}(x -\frac{z}{\aa})-\psi_{j,t}^{-1}(x +\frac{z}{\aa})\big|\cdot \big|\psi_{j,t}^{-1}(x -\frac{z}{\aa}) -\psi_{j,t}^{-1}(x +\frac{z}{\aa})\big|
\end{multline}
Thanks to the mean value theorem (used twice for the first factor and once for the second), we can write:
$$
N_j \leq (\frac{|z|}{\aa})^2 \|D^2 \psi_{j,t}^{-1}\|_{L^\infty} \times(\frac{|z|}{\aa}) \|D \psi_{j,t}^{-1}\|_{L^\infty},
$$
and using the estimates for the flow (see \eqref{estimflot}), we obtain that
$$
N_j \leq e^{CV} (e^{CV}-1) \frac{2^j|z|^3}{\aa^3}.
$$
On the other hand,
$$
D_j=\frac{|z|}{\aa}(Y_++Y_-)\geq 2 e^{-CV(t)}\frac{|z|}{\aa}.
$$
From these estimates, we easily conclude that:
$$
|Y_--Y_+| \leq \frac{\aa}{|z|} \times e^{CV} (e^{CV}-1) \frac{2^j|z|^3}{\aa^3} \times \frac{1}{2} e^{CV}\frac{\aa}{|z|} \leq e^{2CV} (e^{CV}-1) \frac{2^j|z|}{\aa}.
$$
Then, combined with \eqref{estimYYbof}, we obtain:
\begin{equation}
|Y_--Y_+| \leq e^{2CV} (e^{2CV}-1)\min(1,\frac{2^j|z|}{\aa}).
\label{estimRmmoinsRp}
\end{equation}
and then, thanks to \eqref{estimRRbof}, we obtain that:
\begin{equation}
|R_--R_+| \leq C_{d,\eta} e^{C_dV(t)} (e^{2CV}-1)\min(1,\frac{2^j|z|}{\aa})  e^{|z|\left[\eta+(1-\eta)\left(e^{2CV(t)}-1\right)\right]}  (1+|z|^{-\frac{1}{2}}).
\label{estimRR}
\end{equation}
Gathering \eqref{coeffaciles}, \eqref{estimRpm}, \eqref{estim1moinsRpm} and \eqref{estimRR}, we obtain the result from lemma \ref{estimLinf}. $\blacksquare$
\\

\textbf{End of the proof of lemma \ref{lemIIIa}}. Thanks to lemma \ref{estimLinf}, and proposition \ref{propK} (with the same $\eta$), all that remains is to estimate the following function:
\begin{multline}
\phi(z) e^{|z|\left[\eta+(1-\eta)\left(e^{2CV(t)}-1\right)\right]} \leq C_{d,\eta} \frac{e^{-(1-\eta) |z|}}{|z|^{d-\frac{3}{2}}} e^{|z|\left[\eta+(1-\eta)\left(e^{2CV(t)}-1\right)\right]}\\
\leq C_{d,\eta} \frac{1}{|z|^{d-\frac{3}{2}}} e^{|z|\left[-1+2\eta+(1-\eta)\left(e^{2CV(t)}-1\right)\right]}
\end{multline}
Using that $e^{2CV(t)}-1\leq \frac{1}{2}$, we have
$$
-1+2\eta+(1-\eta)\left(e^{2CV(t)}-1\right) \leq -\frac{1}{2}+\frac{3}{2}\eta \leq -\frac{1}{8}
$$
if we choose $\eta=\frac{1}{4}$ (in fact we need $\eta<\frac{1}{3}$). $\blacksquare$

\textbf{Proof of lemma \ref{lemIIIb}:} as $f_j= \ddj f$, we can write that for all $x$ and $z$,
$$
f_j(x-\frac{z}{\aa})+f_j(x+\ee \frac{z}{\aa})-2f_j(x)= \left(\tau_{-\frac{z}{\aa}} \ddj f +\tau_{\frac{z}{\aa}} \ddj f -2\ddj f\right)(x)
$$
and we refer to \cite{CHVP} (proof of lemma $5$) for the following result which is adapted from \cite{Dbook} (theorems $2.36$ and $2.37$):

\begin{lem}
\sl{Under the same assumtptions, there exists a nonnegative summable sequence of summation $1$, that we will also denote by $(c_j)_{j\in \Z}$, such that:
\begin{multline}
 \|\tau_{-\frac{z}{\aa}} \ddj f +\tau_{\frac{z}{\aa}} \ddj f -2\ddj f\|_{L^2}\\
 \leq C_\sigma 2^{-j\sigma} c_j \max(\aa^{-2}, 2^{-2j}) \|\aa^2(\phia *f-f)\|_{\dot{B}_{2,1}^\sigma} \min(1,\frac{ 2^{2j} |z|^2}{\aa^2})
\label{estimdiff2cas}
\end{multline}
\begin{multline}
 \|\tau_{\frac{z}{\aa}} \ddj f -\ddj f\|_{L^2}\\
 \leq C_\sigma 2^{-j\sigma} c_j \max(\aa^{-2}, 2^{-2j}) \|\aa^2(\phia *f-f)\|_{\dot{B}_{2,1}^\sigma} \min(1,\frac{ 2^j |z|}{\aa})
\label{estimdiff2cas2}
\end{multline}
}
\end{lem}
It appears then clearly that both integrals are bounded by:
$$
C_\sigma 2^{-j\sigma} c_j \|\aa^2(\phia *f-f)\|_{\dot{B}_{2,1}^\sigma} \max(\aa^{-2}, 2^{-2j}) \aa^2 \times\mathbb{I} =C_\sigma 2^{-j\sigma} c_j \|\aa^2(\phia *f-f)\|_{\dot{B}_{2,1}^\sigma} \max(1, \frac{\aa^2}{2^{2j}}) \times\mathbb{I},
$$
where we define:
\begin{multline}
\mathbb{I}\overset{def}{=} \int_{\R^d} e^{-\frac{|z|}{8}} \frac{1+|z|^{-\frac{1}{2}}}{|z|^{d-\frac{3}{2}}} \min(1,\frac{ 2^{2j} |z|^2}{\aa^2}) dz =C_d \int_0^\infty e^{-\frac{r}{8}} \frac{1+r^{-\frac{1}{2}}}{r^{d-\frac{3}{2}}} \min(1,\frac{ 2^{2j} r^2}{\aa^2}) r^{d-1}dr\\
=C_d \int_0^\infty e^{-\frac{r}{8}} (1+r^{\frac{1}{2}}) \min(1,\frac{ 2^{2j} r^2}{\aa^2}) dr =C_d(\mathbb{I}_1 +\mathbb{I}_2),
\end{multline}
where we have splitted the integral into:
$$
\mathbb{I}_1= \frac{ 2^{2j}}{\aa^2} \int_0^\frac{\aa}{2^j} e^{-\frac{r}{8}} r^2(1+r^{\frac{1}{2}}) dr \quad \mbox{and} \quad \mathbb{I}_2= \int_\frac{\aa}{2^j}^\infty e^{-\frac{r}{8}} (1+r^{\frac{1}{2}}) dr.
$$
As in \cite{CHVP}, the second integral is easily bounded: as there exists a constant $C>0$ such that for all $x>0$, $(1+r^{\frac{1}{2}}) e^{-\frac{r}{8}} \leq C e^{-\frac{r}{16}}$, we have
$$
\mathbb{I}_2\leq C \int_\frac{\aa}{2^j}^\infty e^{-\frac{r}{16}} dr \leq 16C  e^{-\frac{1}{16}\frac{\aa}{2^j}}.
$$
For the second integral, if we use the same argument, we get the estimate:
\begin{equation}
\mathbb{I}_1\leq C \frac{2^{2j}}{\aa^2} \int_0^\frac{\aa}{2^j} e^{-\frac{r}{16}} dr \leq 16C \frac{ 2^{2j}}{\aa^2} \min(1,\frac{\aa}{2^j})
\label{estI1bof}
\end{equation}
which, multiplied by $\max(1, \frac{\aa^2}{2^{2j}})$ gives a resulting term in $\max(1, \frac{2^j}{\aa})$ that is not summable for high frequencies. To overcome this difficulty we simply write that there is a constant $C'>0$ such that for all $r\in[0,\frac{\aa}{2^j}]$, 
$$
r^2(1+r^{\frac{1}{2}}) e^{-\frac{r}{8}} \leq \frac{\aa}{2^j} r(1+r^{\frac{1}{2}}) e^{-\frac{r}{8}} \leq C' \frac{\aa}{2^j} e^{-\frac{r}{16}}.
$$
and then we obtain:
$$
\mathbb{I}_1\leq C' \frac{2^{2j}}{\aa^2} \frac{\aa}{2^j} \int_0^\frac{\aa}{2^j} e^{-\frac{r}{16}} dr,
$$
which is obviously interesting only when $\frac{\aa}{2^j}\leq 1$ that is for high frequencies: in this case it cancels the diverging $\max(1, \frac{2^j}{\aa})$. But in low frequencies this new term is much bigger that the one from \eqref{estI1bof}, so that we have to combine both estimates and finally get that:
$$
\mathbb{I}_1\leq C \frac{2^{2j}}{\aa^2} \min(1,\frac{\aa}{2^j}) \int_0^\frac{\aa}{2^j} e^{-\frac{r}{16}} dr\leq 16C \frac{2^{2j}}{\aa^2} \min(1,\frac{\aa^2}{2^{2j}}).
$$
We are now able to estimate $\mathbb{I}$:
$$
\mathbb{I}\leq C_d\left(\frac{2^{2j}}{\aa^2} \min(1,\frac{\aa^2}{2^{2j}})+ e^{-\frac{1}{16}\frac{\aa}{2^j}}\right),
$$
and then
\begin{multline}
C_\sigma 2^{-j\sigma} c_j \|\aa^2(\phia *f-f)\|_{\dot{B}_{2,1}^\sigma} \max(1, \frac{\aa^2}{2^{2j}}) \times\mathbb{I} \leq \\
C_\sigma 2^{-j\sigma} c_j \|\aa^2(\phia *f-f)\|_{\dot{B}_{2,1}^\sigma} \left(\max(1, \frac{\aa^2}{2^{2j}}) \frac{2^{2j}}{\aa^2} \min(1,\frac{\aa^2}{2^{2j}})+ \max(e^{-\frac{1}{16}\frac{\aa}{2^j}}, \frac{\aa^2}{2^{2j}}e^{-\frac{1}{16}\frac{\aa}{2^j}})\right)\\
\leq C_\sigma 2^{-j\sigma} c_j \|\aa^2(\phia *f-f)\|_{\dot{B}_{2,1}^\sigma} ,
\end{multline}
which concludes the proof of lemma \ref{lemIIIb}. $\blacksquare$

\subsubsection{End of the proof of theorem \ref{apriori}}

As this part is strictly the same as in \cite{CHVP} (section $3.5$) we will not give details (in particular we refer to \cite{CD} or \cite{CHVP} for estimates on $R_j^{1(2,3)}$ and to \cite{Vishik}). Going back to system \eqref{systchanged}, we use the linear estimates from proposition \ref{estimlinloc}: for all $l\in \Z$, as $\tilde q_j (0)= q_j (0, \psi_{j,0}(.))=q_j(0)$,
\begin{multline}
 \|\ddl \tilde u_j\|_{L_t^\infty L^2}+\nu_0 2^{2l} \|\ddl \tilde u_j\|_{L_t^1 L^2}+(1+\nu 2^l)\left(\|\ddl \tilde q_j\|_{L_t^\infty L^2} +\nu\min(\frac{1}{\ee^2}, 2^{2l})\|\ddj \tilde q_j\|_{L_t^1 L^2}\right)\\
\leq C_{p, \frac{\nu^2}{4\kappa}} \Bigg[ (1+\nu 2^l)\|\ddl q_j (0)\|_{L^2} +\|\ddl u_j(0)\|_{L^2} + (1+\nu 2^l)\|\ddl \tilde f_j+ \ddl R_j^1\|_{L_t^1 L^2}\\
+\|\ddl \tilde g_j +\ddl R_j^2 +\ddl R_j^3 +\kappa \ddl R_j\|_{L_t^1 L^2}\Bigg].
\label{estimlinloc2}
\end{multline}
Thanks to \eqref{energieBF2}, all we need is to estimate the high frequencies, that is $(q_j, u_j)$ for $j\geq 0$. For this, we define some $N_0\in \Z$ (that will be fixed later), and write:
$$
\|q_j\|_{L^2} =\|\tilde q_j \circ \psi_{j,t}^{-1}\|_{L^2} \leq e^{CV} \|\tilde q_j\|_{L^2} \leq e^{CV} \left( \|\dot{S}_{j-N_0}\tilde q_j\|_{L^2} +\sum_{l\geq j-N_0} \|\ddl \tilde q_j\|_{L^2}\right).
$$
We refer to \cite{Dbook} lemma $2.6$, or \cite{Dlagrangien}, lemma $A.1$) for the following classical estimates:
$$
\|\dot{S}_{j-N_0}\tilde q_j\|_{L^2}\leq C e^{CV} \left(e^{CV}-1+ 2^{-N_0} e^{CV}\right) \|q_j\|_{L^2},
$$
so that going back to $U_j$ (we refer to \eqref{Uj} for the definition), we can write that for all $j\geq 0$,
\begin{multline}
 U_j(t) \leq C e^{CV}\Bigg[\left(e^{CV}-1+ 2^{-N_0} e^{CV}\right) U_j(t) + \max(1,2^{2N_0}) C_{p, \frac{\nu^2}{4\kappa}}\times\\
\sum_{l\geq j-N_0} \Big((1+\nu 2^l)\|\ddl q_j (0)\|_{L^2} +\|\ddl u_j(0)\|_{L^2} + (1+\nu 2^l)\Big(\|\ddl \tilde f_j\|_{L_t^1 L^2}+ \|\ddl R_j^1\|_{L_t^1 L^2}\Big)\\
+\|\ddl \tilde g_j\|_{L_t^1 L^2} +\|\ddl R_j^2\|_{L_t^1 L^2} +\|\ddl R_j^3\|_{L_t^1 L^2} +\kappa \|\ddl R_j\|_{L_t^1 L^2}\big)\Bigg],
\label{estimlinloc3}
\end{multline}
We refer to \cite{CHVP} for the estimates on the remainder terms $R_j^{1(2,3)}$. Recall that the previous section is needed to estimate $R_j$. If we fix $N_0>0$ large enough and $t$ so that:
\begin{equation}
 \begin{cases}
\vspace{0.2cm}
\displaystyle{\frac{9}{4} C\cdot 2^{-N_0} \leq \frac{1}{8},}\\
\vspace{0.2cm}
\displaystyle{\frac{9}{4}C(e^{CV}-1)\left(1+2^{5N_0} C_{p, \frac{\nu^2}{4\kappa}} ( \frac{1+|\lambda+\mu|+\mu+\nu}{\nu_0} +\frac{1}{\nu^2})\right) \leq \frac{1}{8}}
\end{cases}
\label{Cond2}
\end{equation}
we obtain that for all $j\geq 0$,
\begin{multline}
 U_j(t) \leq 3 C 2^{5N_0} C_{p, \frac{\nu^2}{4\kappa}} \Bigg((1+\nu 2^j)\|q_j (0)\|_{L^2} +\|u_j(0)\|_{L^2}\\
+(1+\nu 2^j)\|F_j\|_{L_t^1 L^2} +\|G_j\|_{L_t^1 L^2} +\int_0^t 2^{-j(s-1)} c_j(\tau) \|\nabla v(\tau)\|_{\dot{B}_{2,1}^\fd} U(\tau)d\tau\\
+\frac{\kappa}{\nu^2} (V+e^{2CV}-1) \int_0^t c_j(\tau) 2^{-j(s-1)} \nu^2\|\frac{\phi_\ee*\nabla q-\nabla q}{\ee^2}\|_{\dot{B}_{2,1}^{s-1}}d\tau \Bigg).
\end{multline}
Now, if $t$ is so small that:
\begin{equation}
 3 C 2^{5N_0} C_{p, \frac{\nu^2}{4\kappa}} (V+e^{2CV}-1) \leq \frac{1}{2} \frac{\nu^2}{\kappa},
\label{Cond3}
\end{equation}
and if we take $K=(2C_{p,\frac{\nu^2}{4\kappa}})^{-1}$ in \eqref{energieBF2}, then sum over $j\in \Z$, we end up with
\begin{multline}
 U(t) \leq \frac{U(t)}{2} +C_{p,\frac{\nu^2}{4\kappa}} \Bigg(U(0)+ \|F\|_{L_t^1 \dot{B}_{2,1}^{s-1}} +\nu\|F\|_{L_t^1 \dot{B}_{2,1}^s} +\|G\|_{L_t^1 \dot{B}_{2,1}^{s-1}}\\
+(\frac{1+|\lambda+\mu|+\mu+\nu}{\nu_0}+\max(1, \frac{1}{\nu^3})) \int_0^t W'(\tau) U(\tau)\Bigg)
\end{multline}
where
\begin{equation}
 V(t)\overset{def}{=}\int_0^t \|\nabla v(\tau)\|_{L^\infty}d\tau \leq W(t)\overset{def}{=}\int_0^t (\|\nabla v(\tau)\|_{\dot{B_{2,1}^\fd}} +\|v(\tau)\|_{\dot{B_{2,1}^\fd}}^2) d\tau. 
\label{defW}
\end{equation}
and thanks to the Gronwall lemma, we obtain that for $t$ small enough (satisfying conditions \eqref{Cond1}, \eqref{Cond2} and \eqref{Cond3}),
\begin{multline}
 U(t) \leq 2C_{p,\frac{\nu^2}{4\kappa}} \left(U(0)+ \|F\|_{L_t^1 \dot{B}_{2,1}^{s-1}} +\nu\|F\|_{L_t^1 \dot{B}_{2,1}^s} +\|G\|_{L_t^1 \dot{B}_{2,1}^{s-1}}\right)\\
\times e^{\displaystyle{2 C_{p,\frac{\nu^2}{4\kappa}} (\frac{1+|\lambda+\mu|+\mu+\nu}{\nu_0}+\max(1, \frac{1}{\nu^3}))W(t)}}
\end{multline}
Then we globalize the result as in \cite{CHVP}. $\blacksquare$

\subsection{Proof of theorem \ref{thexist}}

\subsection{Existence and uniqueness}

As explained in the introduction, once we have defined the interaction potential $\phia$, we can follow the very same methods as in \cite{CH}: using energy methods gives apriori estimates on the advected linear system. then the proof for existence and uniqueness is classical and follows the lines of \cite{Dbook} (section 10.2.3) for the compressible Navier-Stokes system (see also \cite{DD}, \cite{Dinv}). In order to use the classical Friedrichs approximation, we define the frequency truncation operator $J_n$ by: for all $n\in \N$ and for all $g\in L^2(\R^d)$,
$$
J_n g= \mathcal{F}^{-1}\left(\textbf{1}_{2^{-n}\leq |\xi|\leq C_0 2^n}(\xi) \hat{g}(\xi)\right),
$$
and the following approximated system (we omit the dependency in $\aa$ for more simplicity):
$$
\begin{cases}
\begin{aligned}
&\d_t q_n+ J_n\left(J_n u_n.\nabla J_n q_n\right)+ J_n\div u_n= F_n,\\
&\d_t u_n+ J_n\left(J_n u_n.\nabla J_n u_n\right) -\cA J_n u_n+P'(\overline{\rho}).\nabla J_n q_n-\kappa \aa^2\nabla(\phia*J_n q_n-J_n q_n)= G_n,\\
\end{aligned}
\end{cases}
$$
where
$$
\begin{cases}
 F_n=-J_n\big(J_n q_n.\div J_n u_n\big)\\
G_n= J_n\big(K(J_n q_n).\nabla J_n q_n-I(J_n q_n) \cA J_n u_n\big)
\end{cases}
$$
It is easy to check that it is an ordinary differential equation in $L_n^2\times(L_n^2)^d$, where $L_n^2=\{u\in L^2(\R^d), \quad J_n u=u\}$. Getting uniform estimates implies global lifespan if the initial data is small enough. Then classical compactness arguments give existence of a global solution.

\begin{rem}
\sl{Let us emphasize that the precise estimates proven in the present paper cannot be used in the proof of the existence: indeed the term $J_n\left(u_n.\nabla q_n\right)$ is an obtacle for the lagrangian method to give bounded constants (with respect to $n$). So for the existence we simply use the rough apriori estimates given by the classical energy methods (as they rely on inner products in $L^2$, here $J_n$ has no effect).
}
\end{rem}

To obtain uniqueness for a fixed $\aa$, using \eqref{estimapriori} the computations are close to those for the compressible Navier-Stokes system. As in \cite{Dbook} we need to separate the cases $d\geq 3$ and $d=2$ (the case $d=2$ is more difficult because of endpoints for the remainder estimates in the Littlewood-Paley paradecomposition).

\subsubsection{Uniform estimates}

Using estimate (\ref{estimapriori}) with $s=\fd$, we obtain that for all $t\in \R_+$ (we refer to \cite{CH} for the expressions of the external force terms $F$ and $G$):
\begin{multline}
g(t)=\|(\qa,\ua)\|_{E_\aa^\fd(t)} \leq C_{p,\frac{\nu^2}{4\kappa}} e^{\displaystyle{C_{p,\frac{\nu^2}{4\kappa}} C_{visc}\int_0^t (\|\nabla \ua(\tau)\|_{\dot{B}_{2,1}^\fd}+ \|\ua(\tau)\|_{\dot{B}_{2,1}^\fd}^2)d\tau}}\\
\times\Big(h(0) +\|F\|_{\tilde{L}_t^1 \dot{B}_{2,1}^{\fd-1}}+ \nu\|F\|_{\tilde{L}_t^1 \dot{B}_{2,1}^{\fd}}+ \|G\|_{\tilde{L}_t^1 \dot{B}_{2,1}^{\fd-1}}\Big)
\end{multline}
Let $\eta>0$ be small (it will be fixed later) and assume that:
$$
g(0)\overset{def}{=} \|u(0)\|_{\dot{B}_{2,1}^{\fd-1}}+ \|q(0)\|_{\dot{B}_{2,1}^{\fd-1}} + \nu\|q(0)\|_{\dot{B}_{2,1}^{\fd}} \leq \eta.
$$
Let us now define
$$
T=\sup \{t\in \R_+,\quad g(t)\leq 2 C_{p,\frac{\nu^2}{4\kappa}} g(0)\}.
$$
As $g(0) \leq \eta$, we have $T>0$ ($C>1$) and the aim is to prove by contradiction that $T=\infty$. Assume that $T<\infty$, then we obtain that for all $t\leq T$, denoting by $\C=C_{p,\frac{\nu^2}{4\kappa}}'$ and $C_{visc}'=C_{visc}/\nu_0$
$$
g(t)\leq \C e^{2\C^2 C_{visc}' \eta(1+2\C \eta+e^{2\C^2 C_{visc}' \eta(1+2\C \eta)})} g(0).
$$
So that if $\eta=\eta(p, \frac{\nu^2}{4\kappa}, C_{visc}/\nu_0)>0$ is small enough then for all $t\leq T$, $g(t)< 2C_{p,\frac{\nu^2}{4\kappa}} g(0)$ which contradicts the fact that $T$ is maximal. It then implies that $T=\infty$.

\begin{rem}
\sl{
When for example $\lambda=0$ we have $\mu=\nu=\nu_0$ and for a small $\nu$,  $C_{visc}'\sim \nu^{-4}$, so that the previous condition simply implies that $r=\eta/\nu^4\leq 1$ must be small enough so that $e^{2\C^2 r(1+2\C+ e^{2\C^2})} <2$. A sufficient condition is that $r$ satisfies $e^{6r\C^2 e^{2\C^2}} <2$.
}
\label{condetaexist}
\end{rem}

\subsection{Order parameter estimates}

As already explained (see remark \ref{orderparamexpr}), we have $\ca=\phia * \rho_\aa$, so that $\ca-1=\phia * \qa$ and for all $s\in \R$, $\|\ca-1\|_{\dot{B}_{2,1}^s} \leq \|\qa\|_{\dot{B}_{2,1}^s}$. Moreover:
$$
\|\ca-\rho_\aa\|_{\dot{B}_{2,1}^s} =\|\phia *\qa-\qa\|_{\dot{B}_{2,1}^s} =\frac{1}{\aa^2} \|\qa\|_{\dot{B}_\aa^{s+2,s}} =\Sum_{l\in \Z} 2^{ls} \min(1,\frac{2^{2l}}{\aa^2})\|\ddl \qa\|_{L^2}.
$$
On one hand, thanks to the Lebesgue theorem (for series), for $s\in\{\fd-1, \fd\}$ this involves that:
$$
\|\ca-\rho_\aa\|_{\tilde{L}^\infty(\R_+, \dot{B}_{2,1}^{\fd-1})} +\nu \|\ca-\rho_\aa\|_{\tilde{L}^\infty(\R_+,  \dot{B}_{2,1}^{\fd})} \underset{\aa\rightarrow \infty}{\longrightarrow} 0,
$$
and on the other hand, thanks again to the energy estimates, we end up with:
$$
\nu \|\ca-\rho_\aa\|_{L^1(\R_+, \dot{B}_{2,1}^{\fd-1})} +\nu^2 \|\ca-\rho_\aa\|_{L^1(\R_+, \dot{B}_{2,1}^{\fd})} \leq \frac{C^0}{\aa^2}.
$$
This concludes the proof of theorem \ref{thexist}. $\blacksquare$

\section{Rate of convergence (Theorem \ref{thcv})}

In this section we prove that the solution of $(NSOP_\aa)$ goes to the solution of $(K)$, and we give estimates of the rate of convergence as $\aa$ goes to infinity. Here we follow what we did in \cite{CH}: everything relies on Theorem \ref{apriori}. As already explained, if the initial data satisfy
$$
\|q_0\|_{\dot{B}_{2,1}^{\fd-1}} +\nu\|q_0\|_{\dot{B}_{2,1}^{\fd}} +\|u_0\|_{\dot{B}_{2,1}^{\fd-1}}\leq \eta \leq \min(\eta_K, \eta_{OP}),
$$
then systems $(K)$ and $(NSOP_\aa)$ both have global solutions $(q,u)$ and $(\qa,\ua)$, and with the same notations as before, denoting $\C=C_{p,\frac{\nu^2}{4\kappa}}$ then for all $t\in \R$ we have,
\begin{multline}
g_\aa^\fd(t)\overset{def}{=}\|(\qa,\ua)\|_{E_\aa^\fd(t)} \overset{def}{=}  \|u\|_{\tilde{L}_t^{\infty} \dot{B}_{2,1}^{\fd-1}}+ \|q\|_{\tilde{L}_t^{\infty} \dot{B}_{2,1}^{\fd-1}}+ \nu\|q\|_{\tilde{L}_t^{\infty} \dot{B}_{2,1}^{\fd}}\\
+ \nu_0\|u\|_{\tilde{L}_t^1 \dot{B}_{2,1}^{\fd+1}}+ \nu\|q\|_{\tilde{L}_t^1 \dot{B}_{\aa}^{\fd+1,\fd-1}}+ \nu^2\|q\|_{\tilde{L}_t^1 \dot{B}_{\aa}^{\fd+2,\fd}}\leq \C\eta,
\label{estimop}
\end{multline}
and
\begin{multline}
g^\fd(t)\overset{def}{=}\|u\|_{\tilde{L}_t^{\infty} \dot{B}_{2,1}^{\fd-1}}+ \|q\|_{\tilde{L}_t^{\infty} \dot{B}_{2,1}^{\fd-1}}+ \nu\|q\|_{\tilde{L}_t^{\infty} \dot{B}_{2,1}^{\fd}}+\nu_0\|u\|_{L_t^1 \dot{B}_{2,1}^{\fd+1}}\\
+\nu\|q\|_{L_t^1 \dot{B}_{2,1}^{\fd+1}}+ \nu^2\|q\|_{L_t^1 \dot{B}_{2,1}^{\fd+2}}\leq \C\eta.
\label{estimk}
\end{multline}
As in \cite{CH}, up to an additional forcing term, let us rewrite system $(K)$ with a capillary term as in system $(NSO_\aa)$:
$$
\begin{cases}
\begin{aligned}
&\d_t q+ u.\nabla q+ (1+q)\div u=0,\\
&\d_t u+ u.\nabla u -\cA u+P'(1).\nabla q-\kappa \aa^2 \nabla(\phia*q-q)= K(q).\nabla q- I(q) \cA u+R_\aa,\\
\end{aligned}
\end{cases}
\leqno{(K)}
$$
where the remainder $R_\aa\overset{def}{=}\kappa \nabla\left(\Delta q-\aa^2(\phia*q-q)\right)$ and $K$ and $I$ are defined in the introduction. Let us now write the system satisfied by the difference $(\dq,\du)=(\qa-q, \ua-u)$:
\begin{equation}
 \begin{cases}
\begin{aligned}
&\d_t \dq+ \ua.\nabla \dq+ \div \du=\delta F,\\
&\d_t \du+ \ua.\nabla \du -\cA \du+P'(1).\nabla \dq-\kappa \aa^2\nabla(\phia*\dq-\dq)=\delta G-R_\aa,\\
\end{aligned}
\end{cases}
\end{equation}
where
$$
\begin{cases}
\delta F\overset{def}{=}\sum_{i=1}^3 \delta F_i\\
\delta G\overset{def}{=}\sum_{i=1}^5 \delta G_i
\end{cases}
\mbox{with }
\begin{cases}
 \delta F_1= -\du. \nabla q\\
\delta F_2= -\dq .\div \ua\\
\delta F_3= -q.\div \du
\end{cases}
\mbox{and }
\begin{cases}
 \delta G_1= -\du. \nabla u\\
 \delta G_2= \left(K(\qa)-K(q)\right).\nabla \qa\\
 \delta G_3= K(q).\nabla \dq\\
 \delta G_4= \left(I(\qa)-I(q)\right)\cA \ua\\
 \delta G_5= -I(q)\cA \du.\\
\end{cases}
$$
Except $R_\aa$, all these additionnal terms are exactly the same as in \cite{CH}, so we refer to this article for details on their estimates and we will only focus on what changes. As $\eta$ is small, we can additionnally assume $\eta\leq 1$, let us denote once again by $\C$ a constant only depending on $\frac{\nu^2}{4\kappa}$ and $d$ (that may change from line to line). If we introduce for $h\in[0,1[$
$$
f_\aa(t)\overset{def}{=}\|(\dq,\du)\|_{E_\aa^{\fd-h}(t)},
$$
%
then using \eqref{estimapriori}, we obtain that for all $t\in\R$ (see \cite{CH} for details),
\begin{multline}
f_\aa(t)\leq \C e^{2\C \frac{C_{visc}}{\nu_0}\eta} \\
\times \left[ \int_0^t f_\aa(\tau) F(\tau) d\tau+\int_0^t \|\du\|_{\dot{B}_{2,1}^{\fd-h+1}} \left(\|q\|_{\dot{B}_{2,1}^{\fd-1}}+(1+\nu)\|q\|_{\dot{B}_{2,1}^{\fd}}\right) d\tau +\|R_\aa\|_{L_t^1 \dot{B}_{2,1}^{\fd-h-1}} \right]\\
\leq \C e^{2\C \frac{C_{visc}}{\nu_0}\eta} \left[\int_0^t f_\aa(\tau) F(\tau) d\tau +\eta \frac{1}{\nu_0}(1+\frac{1}{\nu}) f_\aa(t) +\|R_\aa\|_{L_t^1 \dot{B}_{2,1}^{\fd-h-1}}\right].
\end{multline}
where
$$
F(t)=\|q\|_{\dot{B}_{2,1}^{\fd+1}} +\nu\|q\|_{\dot{B}_{2,1}^{\fd+2}} +\|u\|_{\dot{B}_{2,1}^{\fd+1}}+(1+\frac{1}{\nu}) \|\ua\|_{\dot{B}_{2,1}^{\fd+1}}.
$$
\begin{rem}
\sl{We recall that, due to endpoints in the Littlewood-Paley remainder term estimates, as in \cite{CH} we have the condition $h<d-1$, which is why we impose $h<1$ in order to work for any dimension $d\geq 2$.}
\end{rem}
If $\eta$ is so small that $\eta \C e^{2\C \frac{C_{visc}}{\nu_0}\eta} \max(1,\frac{1}{\nu_0})(1+\frac{1}{\nu}) \leq \frac{1}{2}$, then we obtain thanks to the Gronwall lemma,
$$
f_\aa(t)\leq \C e^{2\C \frac{C_{visc}}{\nu_0}\eta} \|R_\aa\|_{L_t^1 \dot{B}_{2,1}^{\fd-h-1}} e^{\C e^{2\C \frac{C_{visc}}{\nu_0}\eta} \int_0^t F(\tau) d\tau}.
$$
Thanks to \eqref{estimop} and \eqref{estimk} we have,
$$
\int_0^t F(\tau) d\tau \leq \left(\frac{1}{\nu}+\frac{1}{\nu_0}(2+\frac{1}{\nu})\right) \C \eta
$$
so that using the condition on $\eta$, we end up with:
\begin{equation}
f_\aa(t)\leq \C (\frac{C_{visc}}{\nu_0}, \frac{\nu^2}{4\kappa},d) \|R_\aa\|_{L_t^1 \dot{B}_{2,1}^{\fd-h-1}}.
\label{estimfalpha}
\end{equation}
\begin{rem}
\sl{For small viscosities, in the case $\lambda=0$, the previous condition on $\eta$ is roughly $\C \frac{\eta}{\nu^2} e^{\C \frac{\eta}{\nu^4}} \leq \frac{1}{2}$ which is obviously implied by the condition required in the existence result (see remark \ref{condetaexist}).
}
\end{rem}
To estimate the remainder in the case $h\in]0,1[$, we simply write that:
$$
\widehat{R_\aa}(\xi)=-i \kappa \xi \frac{|\xi|^4}{\aa^2+|\xi|^2} \widehat{q}(\xi),
$$
which allows us to write:
$$
\|R_\aa\|_{\dot{B}_{2,1}^{\fd-h-1}}=\kappa \Sum_{j\in \Z} 2^{j(\fd-h-1)}\frac{2^{5j}}{\aa^2+2^{2j}}\|q_j\|_{L^2} =\kappa \Sum_{j\in \Z} 2^{j(\fd+2)}\|q_j\|_{L^2} \frac{2^{j(2-h)}}{\aa^2+2^{2j}}
$$
We now consider two cases:
\begin{itemize}
\item If $2^j\geq \aa$ then
$$
\frac{2^{j(2-h)}}{\aa^2+2^{2j}}= 2^{-jh} \frac{2^{2j}}{\aa^2+2^{2j}} \leq \aa^{-h},
$$
\item If $2^j\leq \aa$ then as $h\in [0,1[$
$$
\frac{2^{j(2-h)}}{\aa^2+2^{2j}}= \frac{2^{j(2-a)}}{\aa^2} \leq \frac{\aa^{(2-a)}}{\aa^2} \leq \aa^{-h}.
$$
\end{itemize}
We conclude that:
$$
\|R_\aa\|_{\dot{B}_{2,1}^{\fd-h-1}} \leq \aa^{-h} \|q\|_{\dot{B}_{2,1}^{\fd+2}},
$$
and then
$$
f_\aa(t)\leq \C (\frac{C_{visc}}{\nu_0}, \frac{\nu^2}{4\kappa},d) \|q\|_{L_t^1 \dot{B}_{2,1}^{\fd+2}} \aa^{-h} \leq \C (\frac{C_{visc}}{\nu_0}, \frac{\nu^2}{4\kappa},d) \aa^{-h}.
$$
To complete the proof of the theorem, all that remains is to estimate
$$
\ca-\rho =\phia *\rho_\aa -\rho =\phia *(\rho_\aa-\rho)+ \phia *\rho -\rho =\phia *(\qa-q)+ \phia *q -q
$$
Thanks to what precedes only the last term has to be estimated, and similarly we have
$$
\|\phia *q -q\|_{\dot{B}_{2,1}^{\fd-h-1}}=\Sum_{j\in \Z} 2^{j(\fd-h-1)}\frac{2^{2j}}{\aa^2+2^{2j}}\|q_j\|_{L^2} =\Sum_{j\in \Z} 2^{j(\fd-1)}\|q_j\|_{L^2} \frac{2^{j(2-h)}}{\aa^2+2^{2j}} \leq \aa^{-h} \|q\|_{\dot{B}_{2,1}^{\fd-1}}
$$
and 
$$
\|\phia *q -q\|_{\dot{B}_{2,1}^{\fd-h}} \leq \aa^{-h} \|q\|_{\dot{B}_{2,1}^{\fd}}
$$
For the hybrid norms, we do the same :
\begin{multline}
\|\phia *q -q\|_{\dot{B}_{\aa}^{\fd-h+1, \fd-h-1}} \leq\Sum_{j\in \Z} 2^{j(\fd-h-1)} \min(\aa^2, 2^{2j}) \min(1, \frac{2^{2j}}{\aa^2}) \|q_j\|_{L^2}\\
=\Sum_{j\in \Z} 2^{j(\fd+1)} \|q_j\|_{L^2} \aa^2 2^{-j(h+2)} \min(1, \frac{2^{2j}}{\aa^2})^2 \leq \aa^{-h} \|q\|_{\dot{B}_{2,1}^{\fd+1}}
\end{multline}
so that we finally obtain that if $h\in ]0,1[$:
$$
\|(\rho_\aa-\rho, \ca-\rho, \ua-u)\|_{F_\aa^{\fd-h}} \leq \C (\frac{C_{visc}}{\nu_0}, \frac{\nu^2}{4\kappa},d) \aa^{-h}.
$$
Coming back to \eqref{estimfalpha} in the case $h=0$, we simply write that:
$$
\|R_\aa\|_{\dot{B}_{2,1}^{s}}\leq \kappa \Sum_{j\in \Z} 2^{js}2^{3j} \min(1,\frac{2^{2j}}{\aa^2}) \|q_j\|_{L^2},
$$
and thanks to the Lebesgue theorem for series we obtain that the norm $\|(\rho_\aa-\rho, \ca-\rho, \ua-u)\|_{F_\aa^{\fd-h}}$ goes to zero when $\aa$ goes to infinity, which ends the proof of the theorem. $\blacksquare$

\section{Appendix}

The first part is devoted to a quick presentation of the Littlewood-Paley theory and specific properties for hybrid Besov norms used in this paper. The second section to general considerations on flows.

\subsection{Besov spaces}

\subsubsection{Littlewood-Paley theory}

As usual, the Fourier transform of $u$ with respect to the space variable will be denoted by $\mathcal{F}(u)$ or $\hat{u}$. 
In this section we will briefly state (as in \cite{CD}) classical definitions and properties concerning the homogeneous dyadic decomposition with respect to the Fourier variable. We will recall some classical results and we refer to \cite{Dbook} (Chapter 2) for proofs (and more general properties).

To build the Littlewood-Paley decomposition, we need to fix a smooth radial function $\chi$ supported in (for example) the ball $B(0,\frac{4}{3})$, equal to 1 in a neighborhood of $B(0,\frac{3}{4})$ and such that $r\mapsto \chi(r.e_r)$ is nonincreasing over $\R_+$. So that if we define $\varphi(\xi)=\chi(\xi/2)-\chi(\xi)$, then $\varphi$ is compactly supported in the annulus $\{\xi\in \R^d, c_0=\frac{3}{4}\leq |\xi|\leq C_0=\frac{8}{3}\}$ and we have that,
\begin{equation}
 \forall \xi\in \R^d\setminus\{0\}, \quad \sum_{l\in\Z} \varphi(2^{-l}\xi)=1.
\label{LPxi}
\end{equation}
Then we can define the \textit{dyadic blocks} $(\ddl)_{l\in \Z}$ by $\ddl:= \varphi(2^{-l}D)$ (that is $\hat{\ddl u}=\varphi(2^{-l}\xi)\hat{u}(\xi)$) so that, formally, we have
\begin{equation}
u=\Sum_l \ddl u
\label{LPsomme} 
\end{equation}
As (\ref{LPxi}) is satisfied for $\xi\neq 0$, the previous formal equality holds true for tempered distributions \textit{modulo polynomials}. A way to avoid working modulo polynomials is to consider the set $\cS_h'$ of tempered distributions $u$ such that
$$
\lim_{l\rightarrow -\infty} \|\dot{S}_l u\|_{L^\infty}=0,
$$
where $\dot{S}_l$ stands for the low frequency cut-off defined by $\dot{S}_l:= \chi(2^{-l}D)$. If $u\in \cS_h'$, (\ref{LPsomme}) is true and we can write that $\dot{S}_l u=\Sum_{k\leq l-1} \ddq u$. We can now define the homogeneous Besov spaces used in this article:
\begin{defi}
\label{LPbesov}
 For $s\in\R$ and  
$1\leq p,r\leq\infty,$ we set
$$
\|u\|_{\dot B^s_{p,r}}:=\bigg(\sum_{l} 2^{rls}
\|\ddl u\|^r_{L^p}\bigg)^{\frac{1}{r}}\ \text{ if }\ r<\infty
\quad\text{and}\quad
\|u\|_{\dot B^s_{p,\infty}}:=\sup_{l} 2^{ls}
\|\ddl u\|_{L^p}.
$$
We then define the space $\dot B^s_{p,r}$ as the subset of  distributions $u\in {\cS}'_h$ such that $\|u\|_{\dot B^s_{p,r}}$ is finite.
\end{defi}
Once more, we refer to \cite{Dbook} (chapter $2$) for properties of the inhomogeneous and homogeneous Besov spaces.

In this paper, we mainly work with functions or distributions depending on both the time variable $t$ and the space variable $x.$ We denote by $\cC(I;X)$ the set of continuous functions on $I$ with values in $X.$ For $p\in[1,\infty]$, the notation $L^p(I;X)$ stands for the set of measurable functions on  $I$ with values in $X$ such that $t\mapsto \|f(t)\|_X$ belongs to $L^p(I)$.

In the case where $I=[0,T],$  the space $L^p([0,T];X)$ (resp. $\cC([0,T];X)$) will also be denoted by $L_T^p X$ (resp. $\cC_T X$). Finally, if $I=\R^+$ we alternately use the notation $L^p X.$

The Littlewood-Paley decomposition enables us to work with spectrally localized (hence smooth) functions rather than with rough objects. We naturally obtain bounds for each dyadic block in spaces of type $L^\rho_T L^p.$  Going from those type of bounds to estimates in  $L^\rho_T \dot B^s_{p,r}$ requires to perform a summation in $\ell^r(\Z).$ When doing so however, we \emph{do not} bound the $L^\rho_T \dot B^s_{p,r}$ norm for the time integration has been performed \emph{before} the $\ell^r$ summation.
This leads to the following notation (after J.-Y. Chemin and N. Lerner in \cite{CL}):

\begin{defi}\label{d:espacestilde}
For $T>0,$ $s\in\R$ and  $1\leq r,\rho\leq\infty,$
 we set
$$
\|u\|_{\tilde L_T^\rho \dot B^s_{p,r}}:=
\bigl\Vert2^{js}\|\ddq u\|_{L_T^\rho L^p}\bigr\Vert_{\ell^r(\Z)}.
$$
\end{defi}
It is then possible to define the space $\tilde L^\rho_T \dot B^s_{p,r}$ as the set of  tempered distributions $u$ over $(0,T)\times \R^d$ such that $\lim_{q\rightarrow-\infty}\dot S_q u=0$ in $L^\rho([0,T];L^\infty(\R^d))$ and $\|u\|_{\tilde L_T^\rho \dot B^s_{p,r}}<\infty.$ The letter $T$ is omitted for functions defined over $\R^+.$ 
The spaces $\tilde L^\rho_T \dot B^s_{p,r}$ may be compared with the spaces  $L_T^\rho \dot B^s_{p,r}$ through the Minkowski inequality: we have
$$
\|u\|_{\tilde L_T^\rho \dot B^s_{p,r}}
\leq\|u\|_{L_T^\rho \dot B^s_{p,r}}\ \text{ if }\ r\geq\rho\quad\hbox{and}\quad
\|u\|_{\tilde L_T^\rho \dot B^s_{p,r}}\geq
\|u\|_{L_T^\rho \dot B^s_{p,r}}\ \text{ if }\ r\leq\rho.
$$
All the properties of continuity for the product and composition which are true in Besov spaces remain true in the above  spaces. The time exponent just behaves according to H\"older's inequality. 
\medbreak
Let us now recall a few nonlinear estimates in Besov spaces. Formally, any product of two distributions $u$ and $v$ may be decomposed into 
\begin{equation}\label{eq:bony}
uv=T_uv+T_vu+R(u,v), \mbox{ where}
\end{equation}
$$
T_uv:=\sum_l\dot S_{l-1}u\ddl v,\quad
T_vu:=\sum_l \dot S_{l-1}v\ddl u\ \hbox{ and }\ 
R(u,v):=\sum_l\sum_{|l'-l|\leq1}\ddl u\,\dot\Delta_{l'}v.
$$
The above operator $T$ is called a ``paraproduct'' whereas $R$ is called a ``remainder''. The decomposition \eqref{eq:bony} has been introduced by J.-M. Bony in \cite{BJM}. We refer to \cite{Dbook} for properties, and also to \cite{CH} or \cite{CHVP} for paraproduct and remainder estimates for external force terms.

\subsubsection{Complements for hybrid Besov spaces}

As explained, in the study of the compressible Navier-Stokes system with data in critical spaces, the density fluctuation has two distinct behaviours in low and high frequencies, separated by a frequency threshold. This leads to the notion of hybrid Besov spaces and we refer to R. Danchin in \cite{Dinv} or \cite{Dbook} for general hybrid spaces . In this paper we only will use the following hybrid norms:
\begin{multline}
 \|f\|_{\dot{B}_\aa^{s+2,s}} \overset{def}{=} \sum_{j\in \Z} \min(\aa^2, 2^{2j}) 2^{js} \|\ddj f\|_{L^2}\\
 = \Sum_{j\leq \log_2 \aa} 2^{j(s+2)} \|\ddj f\|_{L^2} +\Sum_{j> \log_2 \aa} \aa^2 2^{js} \|\ddj f\|_{L^2}.
\end{multline}
In this formulation, we obviously remark the threshold frequency $\log_2 \aa$ which separates low (parabolically regularized) and high (damped with coefficient $\aa^2$) frequencies. But as we prove in \eqref{normhybride2}, the frequency transition is in fact continuous and the following equivalent formulations show that these norms are completely tailored to our capillary term:
$$
\|f\|_{\dot{B}_\aa^{s+2,s}} \sim \sum_{j\in \Z} \frac{2^{2j}}{1+\frac{2^{2j}}{\aa^2}}  2^{js} \|\ddj f\|_{L^2} \sim \|\aa^2 (\phi_\ee*f-f)\|_{\dot{B}_{2,1}^s}.
 $$
Exactly as in \cite{CHVP}, the non-local capillary term $\aa^2(\phia*\nabla \qa-\nabla \qa)$ has the same regularity as the capillary term from the local model $\nabla \Delta q$: both of them belong to $L_t^1 (\dot{B}_{2,1}^{\fd-2}\cap \dot{B}_{2,1}^{\fd-1})$.\\

\subsection{Estimates for the flow of a smooth vector-field}
In this section, we recall classical estimates for the flow 
of a smooth vector-field with bounded spatial derivatives. We refer to \cite{Dlagrangien} or \cite{CD} for more details. We also refer to \cite{TH1} for the incompressible Navier-Stokes case.
\begin{prop}
\label{p:flow}
Let $v$ be a smooth globally Lipschitz time dependent  vector-field. Let $W(t) :=\int_0^t\|\nabla v(t')\|_{L^\infty}\,dt'.$  Let $\psi_t$ 
satisfy
$$
\psi_t(x)=x+\int_0^t v(t',\psi_{t'}(x))\,dt'.
$$
Then for all $t\in\R,$ the flow $\psi_t$ is a smooth diffeomorphism over $\R^d$ and one has  if $t\geq0,$
$$
\|D\psi_t^{\pm1}\|_{L^\infty}\leq e^{W(t)},
$$
$$
\|D\psi_t^{\pm1}-I_d\|_{L^\infty}\leq  e^{W(t)}-1,
$$
$$
\|D^2\psi_t^{\pm1}\|_{L^\infty}\leq e^{W(t)} \int_0^t\|D^2v(t')\|_{L^\infty} e^{W(t')} dt',
$$
$$
\|D^3\psi_t^{\pm1}\|_{L^\infty}\leq e^{W(t)} \int_0^t\|D^3v(t')\|_{L^\infty} e^{2W(t')} dt' +3\biggl(e^{V(t)}\int_0^t \|D^2v(t')\|_{L^\infty}e^{W(t')} dt'\biggr)^{2}.
$$
\end{prop}
As in \cite{CD} we also introduce the two-parameter flow $(t,t',x)\mapsto X(t,t',x)$ which is (uniquely) defined by
\begin{equation}
\label{flot2param}
X(t,t',x)=x+\int_{t'}^tv\bigl(t'',X(t'',t',x)\bigr)\,dt''.
\end{equation}
Uniqueness for Ordinary Differential Equations entails that
$$
X(t,t'',X(t'',t',x))=X(t,t',x).
$$ 
Hence 
$\psi_t=X(t,0,\cdot)$ and $\psi_t^{-1}=X(0,t,\cdot).$ 
\begin{prop}
 \sl{Under the previous notations, the jacobian determinant of $X$ satisfies:
\begin{equation}
 det(DX(t,t',x))=e^{\int_{t'}^t (\div v)(\tau, X(\tau,t',x))d\tau},
\end{equation}
and
$$
\begin{cases}
det(D\psi_t(x))=e^{\int_{0}^t (\div v)(\tau,\psi_\tau(x))d\tau},\\
det(D\psi_t^{-1}(x))=e^{-\int_{0}^t (\div v)(\tau, X(\tau,t,x))d\tau} =e^{-\int_{0}^t (\div v)(\tau,\psi_\tau\circ \psi_t^{-1}(x))d\tau}.
\end{cases}
$$
}
\label{detjacobien}
\end{prop}
\textbf{Proof:} differentiating \eqref{flot2param} with respect to $x,$ one gets by virtue of the chain rule,
\begin{equation}
\label{flow4} 
DX(t,t',x)=I_d+\int_{t'}^t  Dv(\tau,X(\tau,t',x))\cdot DX(\tau,t',x)\,d\tau.
\end{equation}
This immediately implies that:
$$
\partial_t (DX)(t,t',x)=Dv(t,X(t,t',x))\cdot DX(t,t',x),
$$
and
$$
\partial_t det(DX(t,t',x))=tr\left(Dv(t,X(t,t',x))\right)\cdot det(DX(t,t',x)),
$$
so that we obtain the result. $\blacksquare$.
\\

\subsection{Bessel functions}

In \eqref{phiexpl} we obtained that the interaction potential can be written the following way: 
$$
\phi(x)=\frac{C_d}{|x|^{\fd-1}} K_{\fd-1}(|x|),
$$
where $K$ is the modified Bessel function of the second kind. In this section we will give specific properties for this Bessel function.
\begin{rem}
\sl{
We refer the reader to \cite{AS},\cite{Bowman}, \cite{Luke} or \cite{Watson} (among a very rich litterature) for a profusion of results and refinements for the Bessel functions. In this paper we will restrict to a very limited number of properties of function $K_\nu$ that will be used in section \ref{sectionRj}.}
\end{rem}
Let us begin with general facts on Bessel functions. If $n$ is an integer, the Bessel function $J_n$ represents the $n$-th Fourier coefficient of the function $\theta\mapsto e^{ix \sin \theta}$:
$$
\forall x\in \R, \quad J_n(x)=\frac{1}{2\pi} \int_0^{2\pi} e^{-i n \theta} e^{ix \sin \theta} d\theta =\frac{1}{\pi} \int_0^\pi \cos( n \theta -x \sin \theta) d\theta.
$$
Alternatively, $J_n(x)$ can be seen as the coefficient of $t^n$ in the development of the function $x\mapsto e^{\frac{x}{2}(t-\frac{1}{t})}$ into powers of $t$.
\begin{rem}
\sl{This function can be extended into $J_\nu(z)$ with $z\in \C$ and for $Re(\nu)>-\frac{1}{2}$ but in this paper we will restrict to real variable and an index $\nu\in \Z$ or $\frac{1}{2}+\Z$.}
\end{rem}
For a general index $\nu$, $J_\nu$ solves the following differential equation:
$$
x^2 y''(x)+ xy'(x)+ (x^2-\nu^2) y(x)=0.
$$
If $\nu\notin \Z$, a basis of the space of solutions of this differential equation is given by $(J_\nu, J_{-\nu})$. If $\nu=n$ is an integer, we have $J_{-n}(x)=(-1)^n(x)$, so that one introduced, as a second element of a basis of solutions, the following function (also called the Bessel function of the second kind) :
$$
Y_n(x)=\lim_{\nu \rightarrow n} \frac{J_\nu(x)-(-1)^n J_{-\nu}(x)}{\nu-n}.
$$
In physics, many functions arise which are similarly constructed from the general Bessel function. Let us now consider the case of the modified Bessel functions of the first and second kind: for a general real index $\nu>-\frac{1}{2}$ and a complex variable $z$ we define
$$
I_\nu (z)=e^{-\frac{1}{2}\nu i \pi} J_{\nu}(i z), \quad \mbox{and}\quad K_\nu(z)=\frac{\pi}{2} \frac{I_{-\nu}(z)-I_\nu(z)}{\sin(\nu \pi)}.
$$
We can prove that $(I_\nu, K_\nu)$ is a basis of the space of solutions of the following differential equation:
$$
x^2 y''(x)+ xy'(x)- (x^2+\nu^2) y(x)=0.
$$
We will now focus on $K_{\nu}$ and recall some important properties (we refer to the cited books):
\begin{prop} [\cite{AS},\cite{Bowman} or \cite{Watson}]
\sl{The function $K_\nu$ satisfies the following properties:
\begin{enumerate}
\item $\displaystyle{K_{\pm \frac{1}{2}}(x)=\sqrt{\frac{\pi}{2}} \frac{e^{-x}}{x^{\frac{1}{2}}}}$.
\item For all $x>0$, $K_\nu(x)>0$.
\item $K_\nu(x)\underset{x\rightarrow 0}{\sim}
\begin{cases}
\vspace{0.2mm}
\displaystyle{-\ln x} & \mbox{ if } \nu=0,\\
\displaystyle{\frac{\Gamma(\nu)}{2} \left(\frac{2}{x}\right)^{\nu}} & \mbox{ if } \nu>0.
\end{cases}
$
\item $\displaystyle{K_\nu(x)\underset{x\rightarrow \infty}{\sim} \sqrt{\frac{\pi}{2}} \frac{e^{-x}}{x^{\frac{1}{2}}}}$.
\item For all $x\in \R_+^*$,
$$
K_\nu(x)=K_{\nu+2}(x)-\frac{2(\nu+1)}{x} K_{\nu+1}(x)=K_{\nu-2}(x)-\frac{2(\nu-1)}{x} K_{\nu-1}(x).
$$
\item $K_\nu$ is decreasing from $\R_+^*$ to $\R_+^*$: for all $x\in \R_+^*$, $K_0'(x)=-K_1(x)$ and
$$
K_\nu'(x)=-\frac{1}{2}(K_{\nu-1}(x)+K_{\nu+1}(x))= -K_{\nu-1}(x)-\frac{\nu}{x} K_\nu(x)= \frac{\nu}{x} K_\nu(x)-K_{\nu+1}(x).
$$
\item For all $x>0$ and $\nu>0$,
$$
K_{\nu}(x)=\int_0^\infty e^{-x ch t} ch(\nu t) dt,
$$
which implies that $K_\nu(x)$ is increasing with respect to $\nu>0$ (see \cite{Reudink}).
\end{enumerate}
}
\label{propgenK}
\end{prop}
We are now able to state the following estimates that will be needed in section \ref{sectionRj}.
\begin{prop}
\sl{Assume that $d\geq 2$.
\begin{enumerate}
\item There exists a constant $C_d>0$ such that for all $x>0$ we have:
$$
K_{\fd-1}(x)\geq
\begin{cases}
\vspace{0.2mm}
\displaystyle{C_d \frac{e^{-x}}{x^{\fd-1}}} & \mbox{if } d>2,\\
\displaystyle{C \frac{e^{-x}}{1+\sqrt{x}}} & \mbox{if } d=2.
\end{cases}
$$
\item For all $\eta\in ]0,1[$, there exists a constant $C_{d,\eta}>0$ such that for all $x>0$ we have:
$$
K_{\fd-1}(x) \leq C_{d,\eta} \frac{e^{-(1-\eta) x}}{x^{\frac{d-1}{2}}}, \quad \mbox{and} \quad |K_{\fd-1}'(x)| \leq C_{d,\eta} \frac{e^{-(1-\eta) x}}{x^{\frac{d+1}{2}}}.
$$
\end{enumerate}
}
\label{propK}
\end{prop}
\textbf{Proof:} For the first point, in the case $\nu=\fd-1>0$ ($d\geq 3$) let us study the following positive function $f:x\in \R_+^*\mapsto f(x)=e^x x^\nu K_\nu(x)$. Thanks to point $3$ from proposition \ref{propgenK}, we have
$$
f(x) \underset{x\rightarrow 0}{\rightarrow} \Gamma(\nu)2^{\nu-1}>0,
$$
and for the limit at infinity,
$$
f(x) \underset{x\rightarrow \infty}{\sim} \sqrt{\frac{\pi}{2}} x^{\nu-\frac{1}{2}} \underset{x\rightarrow \infty}{\rightarrow}
\begin{cases}
\vspace{0.2mm}
\displaystyle{\infty} & \mbox{if } \nu>\frac{1}{2},\\
\displaystyle{\sqrt{\frac{\pi}{2}}} & \mbox{if } \nu=\frac{1}{2},\\
0 & \mbox{if } \nu<\frac{1}{2}.
\end{cases}
$$
Using points $1,5,6$ and $7$ from proposition \ref{propgenK}, we compute the derivative of $f$:
 which is positive :
$$
\forall x>0,\quad f'(x)=e^x x^\nu (K_{\nu}(x)-K_{\nu-1}(x))
\begin{cases}
\vspace{0.2mm}
>0 & \mbox{if }\nu>\frac{1}{2},\\
=0 & \mbox{if }\nu=\frac{1}{2},\\
<0 & \mbox{if }\nu<\frac{1}{2}.
\end{cases}
$$
Thus if $\nu\geq\frac{1}{2}$, for all $x\geq 0$, $f(x)\geq \Gamma(\nu) 2^{\nu-1}$, that is when $d\geq 3$,
$$
K_{\fd-1}(x) \geq \Gamma(\fd -1) 2^{\fd-2} \frac{e^{-x}}{x^{\fd-1}}.
$$
Obviously when $\nu=0$ (that is $d=2$) the previous study is useless and we need to consider $g:x\in \R_+^*\mapsto g(x)=e^x (1+\sqrt{x})K_0(x)$. For all $x>0$, $g(x)>0$ and from the previous results, we have:
$$
\begin{cases}
g(x) \underset{x\rightarrow 0}{\rightarrow} +\infty,\\
g(x) \underset{x\rightarrow \infty}{\rightarrow} \sqrt{\frac{\pi}{2}}. 
\end{cases}
$$
So that either the lower bound $m$ of $g$ satisfies $m\geq \sqrt{\frac{\pi}{2}}$, either we have $m< \sqrt{\frac{\pi}{2}}$ and then $m$ is reached in some segment so that as $g(x)>0$, we have $m>0$. In any case we are sure that $m>0$, which ends the proof of point $1$.

The proof of point $2$ is also very elementary: for a fixed $\eta\in]0,1[$, we introduce the function $h:x\in \R_+^*\mapsto f(x)=e^{(1-\eta)x} x^{\nu+\frac{1}{2}} K_\nu(x)$. For any $\nu\geq 0$, we have:
$$
\begin{cases}
h(x) \underset{x\rightarrow 0}{\rightarrow} 0,\\
h(x) \underset{x\rightarrow \infty}{\rightarrow} 0,
\end{cases}
$$
and $h$ is then a bounded function, which implies the first estimate for $\nu=\fd-1$. For the last point, writing:
$$
K_\nu'(x)= \frac{\nu}{x} K_\nu(x)-K_{\nu+1}(x).
$$
and using what precedes, immediately implies that for all $x>0$,
$$
|K_\nu'(x)|\leq C_{d,\eta} \frac{e^{-(1-\eta) x}}{x^{\nu+\frac{3}{2}}},
$$
and the proof is complete. $\blacksquare$
\\

The author wishes to thank Rapha\"el Danchin, Nicolas Fournier, and Fran\c cois Vigneron for useful discussions.


\begin{thebibliography}{}
\bibitem{AS} M. Abramowitz and I. A. Stegun, Handbook of Mathematical functions with formulas, graphs and mathematical tables, \textit{National Institute of Standards and Technology, Applied Mathematics Series 55}, 1964.
\bibitem{Dbook} H. Bahouri, J.-Y. Chemin, R. Danchin. Fourier analysis and nonlinear partial differential equations, \textit{Grundlehren der mathematischen Wissenschaften}, \textit{343}, \textit{Springer Verlag}, 2011.
\bibitem{BJM} J.-M. Bony, Calcul symbolique et propagation des singularit\'es pour
les \'equations aux d\'eriv\'ees partielles non lin\'eaires, \textit{Annales
Scientifiques de l'\'ecole Normale Sup\'erieure}. 14 (1981)
209-246.
\bibitem{Bowman} F. Bowman, Bessel functions, \textit{Dover Publications, New York,} 1958
\bibitem{BLR} D. Brandon, T. Lin and R. C. Rogers, Phase transitions and hysteresis in local and order-parameter models, \textit{Meccanica}, \textbf{30 (1995)}, 541-565.
\bibitem{CD}
F. Charve, R. Danchin, A global existence result for the compressible Navier-Stokes equations in the critical Lp framework, \textit{Arch. Ration. Mech. Anal.} \textbf{198}(1) (2010), 233-271.
\bibitem{CH}
F. Charve, B. Haspot, Convergence of capillary fluid models: from the non-local to the local Korteweg model, \textit{Indiana U. Math. J.}, \textbf{60}(6), 2011.
\bibitem{CH1}
F. Charve, B. Haspot, Existence of global strong solution and vanishing capillarity-viscosity limit in one dimension for the Korteweg system, \textit{to appear in SIAM J. on Math. Analysis}).
\bibitem{CHVP}
F. Charve, B. Haspot, On a Lagrangian method for the convergence from a non-local to a local Korteweg capillary fluid model, \textit{submitted (in revision)}
\bibitem{CL}
J.-Y. Chemin and N. Lerner, Flot de champs de vecteurs non lipschitziens et \'equations de Navier-Stokes, \textit{J.Differential Equations}, 121 (1992) 314-328.
\bibitem{5CR}
F. Coquel, D. Diehl, C. Merkle and C. Rohde, Sharp and diffuse interface methods for phase transition problems in liquid-vapour flows. Numerical Methods for Hyperbolic and Kinetic Problems, 239-270,  \textit{IRMA Lect. Math. Theor.  Phys,7,Eur. Math. Soc, Z$\ddot{\mbox{u}}$rich}, 2005.
\bibitem{Dinv}
R. Danchin, Global existence in critical spaces for compressible Navier-Stokes equations, \textit{Inventiones Mathematicae}, {\bf 141} (2000), pages 579-614.
\bibitem{Dmach}
R. Danchin, Zero Mach number limit in critical spaces for compressible Navier-Stokes equations, \textit{Annales Scientifiques de l'Ecole Normale Sup\'erieure}, {\bf 35}, pages 27-75 (2002).
\bibitem{DD}
R. Danchin, B. Desjardins Existence of solutions for compressible fluid models of Korteweg type, \textit{Annales de l'IHP, Analyse Non Lin\'eaire}, {\bf 18} (2001), 97-133.
\bibitem{Dlagrangien}
R. Danchin, Uniform estimates for transport-diffusion equations, \textit{J. Hyp. Diff. Eq.}, \textbf{4}(1), 1-17 (2007).
\bibitem{3DS}
J.E. Dunn and J. Serrin, On the thermomechanics of interstitial working, \textit{Arch. Rational Mech. Anal.} {\bf 88(2)} (1985) 95-133.
\bibitem{Gibbs} J.W Gibbs, On the equilibrium of heterogeneous substances, \textit{Transactions of the Connecticut Academy, III.}, Oct. 1875-May 1876, 108-248 and May 1877-July 1878, 343-524.
\bibitem{arma} B. Haspot, Existence of global strong solutions in critical spaces for barotropic viscous fluids,\textit{Arch. Rational Mech. Anal}, {\bf 202}, Issue 2 (2011), Page 427-460.
\bibitem{Has1}
B. Haspot, Cauchy problem for viscous shallow water equations with a term of capillarity, \textit{M3AS}, {\bf 20 (7)} (2010), 1049-1087.
\bibitem{Has2}
B. Haspot, Existence of solutions for compressible fluid models of Korteweg type, \textit{Annales Math\'ematiques Blaise Pascal} {\bf 16}, 431-481 (2009).
\bibitem{Has5} B. Haspot, Cauchy problem for capillarity Van der Waals mode, \textit{Hyperbolic problems: theory, numerics and applications}, 625634, Proc. Sympos. Appl. Math., 67, Part 2, Amer. Math. Soc., Providence, RI, 2009.
\bibitem{TH1}
T. Hmidi, R\'egularit\'e h\"olderienne des poches de tourbillon visqueuses, \textit{Journal de Math\'ematiques pures et appliqu\'ees}, \textbf{84}(11), (2005) 1455-1495.
\bibitem{TH2}
T. Hmidi, S. Keraani, On the global solutions of the super-critical 2D quasi-geostrophic equation in Besov spaces, \textit{Advances in Mathematics} \textbf{214} (2) (2007), 618-638.
\bibitem{TH3}
T. Hmidi, H. Abidi, On the global well-posedness of the critical quasi-geostrophic equation, \textit{SIAM J. Math. Anal.} \textbf{40}(1), 167-185 (2008).
\bibitem{TH4}
T. Hmidi, M. Zerguine, On the global well-posedness of the Euler-Boussinesq system with fractional dissipation, \textit{to appear in Physica D}.
\bibitem{3K}
D.J. Korteweg. Sur la forme que prennent les \'equations du mouvement des fluides si l'on tient compte des forces capillaires par des variations de densit\'e. \textit{Arch. N\'eer. Sci. Exactes S\'er.} II, {\bf 6} :1-24, 1901.
\bibitem{Luke} Y.L. Luke, Integrals of Bessel functions, McGraw-Hill Book company, 1962.
\bibitem{Reudink} D. O. Reudink, On the Signs of the v -Derivatives of the Modified Bessel Functions Iv (x) and Kv (x), \textit{Journal of research of the National Bureau of Standards - B. Mathematical Sciences}, \textbf{72B (4)}, October- December 1968.
\bibitem{5Ro}
C. Rohde, On local and non-local Navier-Stokes-Korteweg systems for liquid-vapour phase transitions. \textit{ZAMM Z. Angew. Math. Mech.} {\bf 85} (2005), no. 12, 839-857.
\bibitem{Rohdehdr}
C. Rohde,   Approximation of Solutions of Conservation Laws by Non-Local Regularization and Discretization, \textit{Habilitation Thesis}, University of Freiburg (2004).
\bibitem{Rohdeorder}
C. Rohde, A local and low-order Navier-Stokes-Korteweg system, \textit{Nonlinear partial differential equations and hyperbolic wave phenomena, 315-337, Contemp. Math.}, \textbf{526 (2010)}, \textit{Amer. Math. Soc., Providence, RI.}
\bibitem{3R}
J.S. Rowlinson, Translation of J.D van der Waals, The thermodynamic theory of capillarity under the hypothesis of a continuous variation of density. \textit{J. Statist. Phys.,} {\bf  20(2)}: 197-244, 1979.
\bibitem{Stein} E. M. Stein, R. Shakarchi, Fourier Analysis, an introduction, \textit{Princeton lectures in Analysis I, Princeton University press}, 2003.
\bibitem{VW}
J.F Van der Waals, Thermodynamische Theorie der Kapillarit\"at unter Voraussetzung stetiger Dichte\"anderung, \textit{Phys. Chem.} {\bf 13}, 657-725 (1894).
\bibitem{Vishik}
M. Vishik, Hydrodynamics in Besov spaces, \textit{Arch. ration. Mech. Anal.}, \textbf{145}, 197-214 (1998)
\bibitem{Watson} Watson, A treatise on the theory of Bessel functions, \textit{Cambridge University Press}, 1922.
\end{thebibliography}
\end{document}